\documentclass[12pt]{article}

\usepackage[top=1in,bottom=1in, left=1in, right=1in]{geometry}

\usepackage{amsmath,amsthm,amssymb,amsopn,amsfonts,pdfpages,url,dsfont}
\usepackage[linesnumbered,ruled,vlined]{algorithm2e}
\usepackage{graphics,relsize}
\usepackage{graphicx}
\usepackage{multirow}
\usepackage{bbm}
\usepackage{authblk}
\usepackage{caption}
\usepackage{subcaption}
\usepackage{setspace} 
\usepackage[colorlinks=true,linkcolor=cyan,citecolor=blue]{hyperref}
\usepackage{natbib}

\newcommand{\twospace}{\renewcommand{\baselinestretch}{1.3}\normalsize}

\newcommand{\p}{\mathbb{P}}
\newcommand{\q}{\vec{\mathbf{1}}}
\newcommand{\e}{\mathbb{E}}
\newcommand{\R}{\mathbb{R}}
\newcommand{\Om}{\Omega}
\newcommand{\OM}{\Lambda}
\newcommand{\vP}{\varPi}
\newcommand{\vD}{{\mathcal D}}
\newcommand{\TT}{\textup{trace}}
\newcommand{\lb}{\underline{\xi}}
\newcommand{\ub}{\overline{\xi}}
\newcommand{\XX}{{\mathcal X}}
\newcommand{\nn}{\mathfrak{n}}

\newcommand{\II}{\bigl [\begin{smallmatrix} I_K & 0 \\ 0 & 0  \end{smallmatrix} \bigr ]}
\newcommand{\EE}{{\mathcal E}}

\newtheorem{theorem}{Theorem}
\newtheorem{lemma}{Lemma}

\theoremstyle{definition}
\twospace

\begin{document}

\title{On Seeded Subgraph-to-Subgraph Matching:\\ \makebox[0.9\textwidth][c]{The ssSGM Algorithm and Matchability Information Theory}}

\author[1]{Lingyao Meng}
\author[2]{Mengqi Lou}
\author[1]{Jianyu Lin}
\author[3]{\\Vince Lyzinski}
\author[1]{Donniell E. Fishkind}
\affil[1]{Johns Hopkins University, Department of Applied Mathematics and Statistics}
\affil[2]{Georgia Institute of Technology, School of Industrial and Systems Engineering}
\affil[3]{University of Maryland, College Park, Department of Mathematics}

\maketitle
\doublespacing
\begin{abstract} The subgraph-subgraph matching problem is, given a pair of graphs and a positive integer~$K$,
to find $K$ vertices in the first graph, $K$ vertices in the second graph, and a bijection between them, so as to
minimize the number of adjacency disagreements across the bijection; it is ``seeded" if some of this bijection is fixed.
The problem is intractable, and we present the ssSGM algorithm, which uses Frank-Wolfe methodology to efficiently
find an approximate solution. Then, in the context of a generalized correlated random Bernoulli graph model,
in which the pair of graphs naturally have a core of $K$ matched pairs of vertices, we provide and prove mild
conditions for the subgraph-subgraph matching problem solution to almost always be the correct $K$ matched pairs of vertices.\\
\noindent
{\bf MSC2020 subject classifications:} 05C60 05C80  90C35

\noindent {\bf Authors' emails:} lmeng2@jh.edu, mlou30@gatech.edu, jlin153@jh.edu, vlyzinsk@umd.edu, and def@jhu.edu
\end{abstract}

\newpage
\section{Background and overview}
\label{sec:intro}

Given a pair of graphs, the graph matching problem seeks to discover a correspondence between the vertex sets of the graphs that minimizes a measure of graph dissimilarity; see, for example, \citet{conte2004thirty,foggia2014graph,gmrev} for surveys of the recent graph matching literature.
This problem of optimally aligning a pair of graphs is practically useful, theoretically challenging, and computationally difficult.
Practically, graph matching methods have been successfully used for aligning protein interaction networks \citet{zaslavskiy2009global,netal}, uncovering bilateral brain structure \citet{pedigo2022bisected}, image alignment in computer vision \citet{zhou2012factorized,pose,sun2020survey}, and de-anonymizing social networks \citet{chiasserini2018anonymizing,patsolic2017vertex}, to name a few applications.
Theoretically, there is a growing literature devoted to uncovering sharp information-theoretic thresholds for graph matchability (i.e., when can graph matching algorithms uncover a latent correspondence across the vertex sets of the graph pair in the presence of noise); see, for example,
\citet{Gross2013PGM,JMLR:v15:lyzinski14a,cullina2016improved,cullina2017exact,wu2021settling,detection, FISHKIND2019295, fishkind2021phantom}.
While much of the matchability work is situated in the setting of Erd\H os-R\'enyi graphs, analogous matchability results have been derived in the stochastic blockmodel setting
\citet{lyzinski2016information,shirani2018matching}; of note in the stochastic blockmodel setting, there is a fascinating interplay between matchability and community recovery in multiple correlated blockmodel graphs, see \citet{racz1,racz2021correlated}.
In parallel to this work on theoretic matchability, there is a branch of work in the literature focused on developing efficient algorithms for matching graphs with theoretical guarantees of asymptotically perfect recovery in the Erd\H os-R\'enyi setting (see, for example, \citet{otter,Grampa,profile}), and the stochastic blockmodel setting (see, for example, \citet{racz2}).

Computationally, when the graphs are allowed to be weighted, directed and loopy, the graph matching problem is equivalent to the NP-hard quadratic assignment problem.
Even the simpler graph isomorphism problem---determining if two graphs are isomorphic---is NP-intermediate (see \citet{babai2016graph} for the remarkable proof of the quasipolynomial time complexity of the problem).
This combination of utility and complexity has led to a huge number of approximate graph matching algorithms in the literature; see, for example, \citet{VogConLyzPodKraHarFisVogPri2015,klau2009new,final,regal,Gross2011OPAN,isorankn,netalign,fast,kergm}.

One particular branch of algorithms (see, for example, \citet{VogConLyzPodKraHarFisVogPri2015,sgm,kergm,saad2021graph})
proceeds by relaxing the combinatorial graph matching problem
to a continuous optimization problem via relaxing the permutation matrix constraint to the set of doubly-stochastic matrices (e.g. see Section~\ref{section:fwm}).
This constrained, continuous optimization problem can then be approximately solved via the Frank-Wolfe method \citet{FW} or a fast approximation thereof (see for example \citet{cuturi2013sinkhorn,kergm,saad2021graph}), and the relaxed solution is projected to the nearest feasible point in the original combinatorial problem.
These Frank-Wolfe-based algorithms empirically provide excellent performance, especially in the presence of seeded vertices (i.e.~when some of a natural correspondence is known a priori).

In the present work, we consider the problem of finding the best-matched pair of 
subgraphs---of a fixed order---in a pair of graphs,  simultaneously optimally aligning these subgraphs, also incorporating any seeds; this {\it seeded subgraph-subgraph matching problem} is described in Section~\ref{sec:ssgmp}. When there is a correlated core in the two graphs to be matched, then the seeded subgraph-subgraph matching problem can be used to identify the core as well as the correspondence between the core vertices across the graphs.

Our setting contrasts with the work in \citet{hu2022graph} that considers the task of finding the best {\it largest} subgraph alignment, hence the number of vertices in the subgraph is not fixed. Their formulation maximizes the number of corresponding adjacencies across two subgraphs as an equivalent surrogate for minimizing the number of adjacency disagreements across two subgraphs; unfortunately, this surrogacy is not equivalent when the number of vertices of the subgraph is not fixed. Indeed, there could be an isomorphism when considering half the vertices of the two graphs, yet three quarters of the vertices of the two graphs might score a higher number of common adjacencies. This is not to be taken as criticism of the approach of  \citet{hu2022graph}, because if
the order of the core is unknown then it is not clear that there is anything better to do. (Indeed, we demonstrated in \citet{rel} that, for graph matching (not subgraph matching) based on solving relaxations, maximizing (a relaxation of) corresponding adjacencies is effective, but minimizing (a relaxation) of adjacency disagreements is not effective.)
In this paper, by contrast, we assume that the number of vertices in the core is known, hence fixed (or, in practice, approximately known).  It should also be noted that the seeding in    \citet{hu2022graph} is unrelated to the seeding in our setting.

{\bf Our contributions are two-fold:} {\bf First}, in Section \ref{sec:ssSGMalg} we develop the ssSGM algorithm to efficiently, approximately solve the seeded subgraph subgraph matching problem; it incorporates several modifications to the Frank-Wolfe methodology of the SGM Algorithm \citet{sgm} for seeded graph matching. Sections~\ref{sec:exp},~\ref{sec:real} illustrate ssSGM with simulations and real-data experiments.
{\bf Our second contribution} is an information-theoretic result in Section~\ref{section:theo};
if there is a latent alignment of a pair of subgraphs across the pair of graphs,
we provide and prove mild conditions where the exact solution of
the subgraph-subgraph matching problem recovers the latent alignment.
These results are of note in that much of the information theoretic thresholds in the literature are developed in the setting of perfectly-overlapping random graphs; 
notable exceptions include
\citet{sussman2018matched,yang2023structural} which are focused on subgraph detection rather than subgraph alignment.
Discussion follows in Section~\ref{sec:discussion}.

\section{Graph Matching and subgraph-subgraph matching}
\label{sec:gmp}

All graphs in this paper are simple (i.e.~no orientation to edges, no loops, no multi-edges). For any positive integer $n$, we use the notation $[n]$ to denote the set of 
integers $\{1,2,3,\ldots,n\}$.

Suppose $G$ and $H$ are graphs with respective vertex sets $V$ and $W$ with the same cardinality, say $n:=|V|=|W|$.
For any bijection $\phi: V \rightarrow W$, a {\it disagreeing pair} is any pair of vertices $v,v' \in V$ such that $v \sim v'$ and $\phi (v) \not \sim \phi (v')$,
or  $v \not \sim v'$ and $\phi (v) \sim \phi (v')$. The {\it graph
matching problem} is to find a bijection $\phi: V \rightarrow W$ with minimum number of disagreeing pairs.

The graph matching problem can be equivalently expressed as follows. 
Say the vertex sets of $G$ and $H$ are respectively $V=\{v_1,v_2,\ldots$ $,v_n\}$ and
$W=\{ w_1,w_2,\ldots,w_n\}$. Let $\vP_n$ denote the set of $n \times n$ permutation matrices. For any real numbers $\delta^{(1)}$, $\delta^{(2)}$, and $\delta^{(3)}$ such that $\delta^{(1)} \ne \delta^{(2)}$, let $A$ and $B$ be  $(\delta^{(1)},\delta^{(2)},\delta^{(3)})$-adjacency matrices for $G$ and $H$ respectively; this means that 
for all $i,j \in [n]$
\begin{eqnarray} \label{eqn:adjacency}
A_{ij}= \left \{
\begin{array}{ll} \delta^{(1)} & \mbox{if } v_i \sim v_j \\
                  \delta^{(2)} & \mbox{if } v_i \not \sim v_j \mbox{ and } i \ne j \\
                  \delta^{(3)} & \mbox{if } i=j
                  \end{array}
\right .,
\end{eqnarray}
and $B_{ij}$ is defined in the same manner, substituting ``$w$" for ``$v$" in the definition of $A_{ij}$.  The graph matching problem can be equivalently expressed as
$\min \| A - XBX^T\|_F^2$ s.t.~$X \in \vP_n$, where $\| \cdot \|_F$ is the Frobenius norm.
The usual choice of $\delta^{(1)}$,
$\delta^{(2)}$, and $\delta^{(3)}$ are,
respectively, $1$, $0$, and $0$, in which case $A$ and $B$ are the standard
``$\{0,1 \}$-adjacency matrices" of $G$ and $H$, respectively. This will turn out to be an unhelpful choice for us here, and we will need to adopt an alternative, as described next in Section \ref{sec:ssgmp}. Specifically, we will take $\delta^{(1)}$,
$\delta^{(2)}$, and $\delta^{(3)}$ as,
respectively, $1$, $-1$, and $-1$; these $(1,-1,-1)$-adjacency matrices will also be called $\pm 1-$adjacency matrices.

The implicit story is that there is a natural correspondence between the two vertex sets that we wish to discover. For example, if the two graphs are modeling two
different social media platforms with the same users in each platform (perhaps the users identify themselves by their real names in one platform and a different identifier in the other platform) and, for any two users, adjacency in a platform
means the respective users communicated in that platform, then the graph matching problem would attempt to cross-identify the users between platforms.

The \textit{seeded} graph matching problem refers to the setting where the natural  correspondence between the vertex sets is partially known a priori, so we restrict the
optimization to bijections that agree with the known correspondences, see
\citet{JMLR:v15:lyzinski14a}.
In the social media story above, this would happen if some users identify themselves by their real name in both platforms. (See \citet{fang2018tractable} for a variant, known as \textit{soft-seeding}, where the a priori known alignment is soft-coded into the search space and is allowed to be updated; this would be useful if the a priori known correspondences were not completely certain. We will not be working with soft-seeding here.)
The practical utility of even a few seeds on matching performance (see \citet{sgm}) has led to seeded graph matching becoming a popular variant of the classical framework; see, for example, \citet{mossel2020seeded,seeded_erkip,shirani2017seeded,sgmnet,sgm_jofc,seeed} for state of the current algorithmic and theoretic seeded graph matching problem literature.

The focus of our work here is a further generalization of this problem setting,
called the subgraph-subgraph matching
problem, which we describe in detail next
in Section \ref{sec:ssgmp}.

\subsection{The seeded subgraph-subgraph matching problem } 
\label{sec:ssgmp}

In this section, we formally present the seeded subgraph-subgraph matching problem,
which is the main focus of our work here.

Suppose $G$ and $H$ are graphs with respective vertex sets $V$ and $W$,
say $m:=|V|$, $n:=|W|$.
Let $K$ be a positive integer such that $K \leq \min \{ m,n \}$. For any two
subsets $\Om \subseteq V$ and $\OM \subseteq W$ such that $| \Om |= | \OM | = K$,
and any bijection $\phi: \Om \rightarrow \OM$, a {\it disagreeing pair} is any
pair of vertices $v,v' \in \Om$ such that $v \sim v'$ and $\phi (v) \not \sim \phi (v')$,
or  $v \not \sim v'$ and $\phi (v) \sim \phi (v')$. The {\it subgraph-subgraph
matching problem} is to find subsets $\Om \subseteq V$ and $\OM \subseteq W$ with $| \Om |= | \OM | = K$,
and to find a bijection $\phi: \Om \rightarrow \OM$, minimizing the number of disagreeing pairs (optimized over all possible
choices of pairs of vertex subsets {\bf and also} all bijections between the corresponding vertex subsets).

Analogously to the previous section, the implicit story is that there is a natural correspondence between a {\it core} of $K$ vertices in each graph,
and we want to estimate the core and the correspondence. For example, if the two graphs are modeling two
different social media platforms, where adjacency of two vertices in a platform
means the respective users communicated in that platform, and if we know that the number of users who are common to
both platforms is $K$ (albeit, each user with different identifiers in the two platforms), then the subgraph-subgraph matching problem would attempt to
find these $K$ common users and to cross-identify them between platforms.

We can equivalently express the problem as follows. Say the vertices are $V=\{v_1,v_2,\ldots$ $,v_m\}$ and
$W=\{ w_1,w_2,\ldots,w_n\}$, and then let $A \in \{1,-1\}^{m \times m}$ and
$B \in \{1,-1 \}^{n \times n}$ be the respective $\pm 1$-adjacency matrices of $G$ and $H$.
(That is, for all indices $i,j$, we have that $A_{ij}$ is $1$ or $-1$ according as $v_i \sim v_j$ or
$v_i \not \sim v_j$, and $B_{ij}$ is $1$ or $-1$ according as $w_i \sim w_j$ or
$w_i \not \sim w_j$.) For any positive integers $c,d,e$ such that $e\leq \min \{c,d\}$, let $\vP_{c,d,e}$ denote the set of $c \times d$ binary matrices with at most one $1$ per row
and at most one $1$ per column, with the total number of $1$'s being $e$; that is,
$\vP_{c,d,e} := \{ X \in \{0,1\}^{c \times d}: X \q_d \leq \q_c, \q_c^TX \leq \q_d^T, \q_c^T X \q_d=e \}$, where
$\q$ denotes a vector of ones with length given by its subscript. The
subgraph-subgraph matching problem can be expressed as
$\min \| A - XBX^T\|_F^2$ s.t.~$X \in \vP_{m,n,K}$; indeed, for all $i\in [m]$ and $j\in [n]$ such that
$X_{ij}=1$, if we understand this to say that $v_i \in \Om$ and $w_j \in \OM$ and $\phi (v_i)=w_j$, then the objective function
is eight times the number of disagreeing pairs plus a constant (this constant is $m^2-K^2$), so we are
essentially minimizing the number of disagreeing~pairs.

Next, we  expand the objective function $\| A - XBX^T\|_F^2=\|A\|^2_F- 2 \TT A^TXBX^T +\|XBX^T\|^2_F$; note that, when restricted to $X \in \vP_{m,n,K}$, we have that   $\|A\|_F^2$ and
$\|XBX^T\|_F^2$ are constants, respectively equal to $m^2$ and $K^2$, and note that $A$ is symmetric, from which we get the
following equivalent formulation of the subgraph-subgraph matching problem:
\begin{eqnarray}
\max_{X \in \vP_{m,n,K}}\TT AXBX^T \ \ .
\end{eqnarray}
{\it Note that if we had instead used standard $\{0,1\}$-adjacency matrices, rather than $\pm 1$-adjacency
matrices, then this expression would {\bf not} be equivalent to the subgraph-subgraph matching problem.}

In the social media story above, sometimes we will know a priori a few members of the core and their correspondence
across the media platforms. With this in mind, the
{\it seeded subgraph-subgraph matching problem} is where the subgraph-subgraph matching optimization is
done when fixing $s$ vertices of $\Om$, $s$ vertices of $\OM$, and $\phi$ is fixed in mapping the former fixed vertices to the latter fixed vertices in some particular manner.
Without loss of generality, the fixed vertices (``seeds") are $v_1,v_2,\ldots, v_s$ and $w_1,w_2,\ldots, w_s$ and, without loss of generality,
for each $i\in [s]$, we fix $\phi (v_i)=w_i$. Thus we can equivalently express:
\begin{eqnarray} \label{sssmp} (\textup{seeded subgraph-subgraph matching problem}) \ \
\max_{X \in \vP_{m-s,n-s,K-s}}\TT A(I_s \oplus X)B(I_s \oplus X)^T \ \ ,
\end{eqnarray}
where ``$\oplus$" denotes the direct sum of matrices (a block-diagonal matrix with the summands as the diagonal blocks).
This problem is computationally intractable;
indeed, the special cases seeded graph matching problem and graph matching problem
are computationally intractable, as elaborated in Section~\ref{sec:intro}.

Thus we seek approximate solutions that are efficiently
computable. In the next section, we adapt the Frank-Wolfe methodology---that was used for
efficient approximate solution of the seeded graph matching problem in \citet{sgm}---to our setting of seeded subgraph-subgraph graph matching.

Much of the literature on subgraph matching is focused on matching $G$ to a subgraph of $H$ (i.e., $K=m<n$) where $H$ is often of much larger order than $G$; see, for example, \citet{billion,yang2023structural,moorman2018filtering,moorman2021subgraph,sussman2018matched}.
Our problem is more akin to the $K$-core alignment problem \citet{cullina2019partial} or partial alignment problem \citet{partial,wu2021settling}, wherein we seek to correctly partially align (often all but a vanishing fraction of) a pair of large networks.
However, in our present framework we do not assume the subgraphs we seek to align have any structural properties a priori (i.e., density, average degree, etc.), and we do not posit a latent (or structural) alignment of the remainder of the graphs.

\section{The ssSGM algorithm \label{sec:ssSGMalg}}

In this section we introduce the ssSGM algorithm, which provides an efficiently computable,
approximate solution to the seeded subgraph-subgraph matching problem, by
generalizing the techniques of \citet{VogConLyzPodKraHarFisVogPri2015,sgm,sussman2018matched}. The algorithm is described as follows.
We begin by relaxing the seeded subgraph-subgraph matching problem, and then we apply Frank-Wolfe methodology to approximately
solve the relaxed problem. Finally, we project the relaxed solution to the nearest feasible point
in the original problem, which provides an approximate solution to the seeded subgraph-subgraph
matching problem. We begin by stating two key technical lemmas.

\subsection{Key technical lemmas}

For any positive integers $c,d,e$ such that $e\leq \min \{c,d\}$, let $\vD_{c,d,e}$
denote the set of $c \times d$ nonnegative, real matrices with row sums at most $1$, column sums at most $1$,
and all entries summing to $e$; that is, $\vD_{c,d,e} :=
\{ X \in \R_{\geq 0}^{c \times d}: X \q_d \leq \q_c, \q_c^TX \leq \q_d^T, \q_c^T X \q_d=e \}$. The following
result is quite useful; it is essentially due to Mendelsohn and Dulmage \citet{MendelsohnDulmage} (as extended in \citet{Gilletal} ), and it is a generalization
of the classical Birkhoff-Von Neumann Theorem.
Let $\mathcal{H}(\cdot)$ denote the convex hull of its argument.

\begin{lemma} \label{lemma:hull}
For any positive integers $c,d,e$ such that $e\leq \min \{c,d\}$, it holds  $\mathcal{H}(\vP_{c,d,e})=\vD_{c,d,e}$.
\end{lemma}

\noindent
Next, let $M$ be any given matrix in $\R^{c \times d}$. The {\it generalized linear assignment problem} is defined as
\begin{eqnarray} \label{glap}
\max_{X \in \vP_{c,d,e}} \TT M^TX
\end{eqnarray}
The special case where $c=d=e$ is the well-known linear assignment problem, and is solvable in time $O(e^3)$ by the
Hungarian Algorithm \citet{hungarian,jonker1987shortest}.
 Note that (\ref{glap}) is equivalent to $\max_{X \in \vD_{c,d,e}}\TT M^TX$; this is an immediate consequence of
Lemma \ref{lemma:hull} and since, for any linear function, the function value of a convex combination
of vectors is the convex combination of the vectors' function values.

Let $\lb$ be any real number less than all of the entries of $M$, and let $\ub$ be any real number greater than
all of the entries of $M$. Define the matrix $\mathcal{M} \in \R^{(c+d-e)\times (c+d-e)}$ such that
$\mathcal{M}= \bigl [\begin{smallmatrix} M & \ub E \\ \ub E & \lb E \end{smallmatrix} \bigr ]$, where $E$ denotes a
matrix of all $1$s, appropriately sized. In \ref{app:pfs} we prove the following lemma.

\begin{lemma} \label{lemma:assign}
For any positive integers $c,d,e$ such that $e\leq \min \{c,d\}$, and for any $M \in \R^{c \times d}$, suppose
$\mathcal{X}^* \in \vP_{c+d-e,c+d-e,c+d-e}$ is a solution of the linear assignment problem $$\max_{\XX \in \vP_{c+d-e,c+d-e,c+d-e}}\TT \mathcal{M}^T \XX.$$
Let $X^*$ be the $c \times d$ upper-left corner submatrix of $\mathcal{X}^*$; i.e.~induced by the first $c$ rows and first $d$
columns of $\mathcal{X}^*$. Then $X^*$ is a solution to the generalized linear assignment problem
$\max_{X \in \vP_{c,d,e}} \TT M^TX$.
\end{lemma}

In the next section, Lemma \ref{lemma:assign} will provide us with a means 
in Algorithm \ref{alg:glap} of solving the generalized linear assignment problem, a critical subroutine of the ssSGM algorithm. 

\subsection{Solving the relaxed problem with Frank-Wolfe methodology \label{section:fwm} }

In order to efficiently obtain an approximate solution to the
seeded subgraph-subgraph matching problem as expressed in (\ref{sssmp}), we begin by expressing its relaxation:
\begin{eqnarray} (\textup{relaxed seeded subgraph-subgraph matching problem}) \ \
\max_{X \in \vD_{m-s,n-s,K-s}}\TT A(I_s \oplus X)B(I_s \oplus X)^T  .    \label{problem:relaxed}
\end{eqnarray}
We will apply Frank-Wolfe  \citet{FW} methodology to approximately solve this relaxed problem; afterwards, we will project the
relaxed solution back  to the original feasible region.

In general, when maximizing any objective function over any polytope, Frank-Wolfe methodology consists of
generating a sequence of iterates, which are points in the polytope, terminating when an iterate does not change much
from the previous iterate. The first iterate is an arbitrary point in the polytope and, given any iterate, the next iterate
is calculated from the current iterate in two steps.
The first step is to linearize the objective function about the current iterate, and to maximize this
linear function over the polytope; this is simply a linear programming problem.
The second step is to maximize the (original) objective function along (only) the line segment
between the current iterate and  the linear program  solution from the first step;
the solution to this one-dimensional optimization problem is the next iterate.

The objective function for the seeded subgraph-subgraph matching problem is the
quadratic function $f(X):=\TT A(I_s \oplus X)B(I_s \oplus X)^T$, and we want to maximize it over the set
$\vD_{m-s,n-s,K-s}$, which is clearly a polytope. A natural choice for the first Frank-Wolfe iterate is the
$(m-s) \times (n-s)$ matrix with all entries $\frac{K-s}{(m-s)(n-s)}$, which is clearly in $\vD_{m-s,n-s,K-s}$.

Given any current Frank-Wolfe iterate, say $Z \in \vD_{m-s,n-s,K-s}$, the linearization of $f$ about $Z$ is given by
$\TT [\nabla f(Z)]^TX$ plus a fixed constant that can be ignored when maximizing. If we partition
$A= \bigl [\begin{smallmatrix} A^{(11)} & A^{(21)^T} \\ A^{(21)} & A^{(22)} \end{smallmatrix} \bigr ]$ and
$B= \bigl [\begin{smallmatrix} B^{(11)} & B^{(12)} \\ B^{(12)^T} & B^{(22)} \end{smallmatrix} \bigr ]$, where
$A^{(21)} \in \{ -1,1 \}^{(m-s) \times s}$, $B^{(12)} \in  \{-1,1 \}^{s \times (n-s)}$,
$A^{(22)} \in \{ -1,1 \}^{(m-s)\times (m-s)}$, and $B^{(22)} \in \{ -1,1 \}^{(n-s)\times (n-s)}$ then
$\nabla f (X)= 2A^{(21)}B^{(12)}+2A^{(22)}XB^{(22)}$, hence the first step to computing the next Frank-Wolfe
iterate is to solve the generalized linear assignment problem
\begin{eqnarray}
\max_{X \in \vD_{m-s,n-s,K-s}}\TT \left ( A^{(21)}B^{(12)}+A^{(22)}ZB^{(22)} \right )^TX,
\end{eqnarray}
which can be done efficiently, in time $O(\max\{m,n\}^3)$, via the Hungarian Algorithm, thanks to Lemma \ref{lemma:assign} (since the matrix 
for the Hungarian Algorithm is of size $O(\max\{m,n\}$-by-$O(\max\{m,n\}$).
Say the solution $X$ to this linear program is $X^* \in \vD_{m-s,n-s,K-s}$. The second step to computing the
next Frank-Wolfe iterate is maximizing $f$ along the line
segment $\{ \lambda X^*+(1-\lambda)Z : \lambda \in [0,1] \}$; this is trivial to do since $f$ is a
quadratic function in $\lambda$ on this line segment. Say the optimal $\lambda$ is $\lambda^*$. The the next Frank-Wolfe
iterate is defined to be $\lambda^* X^*+(1-\lambda^*)Z \in \vD_{m-s,n-s,K-s} $.

When the Frank-Wolfe algorithm terminates, say that $Z \in \vD_{m-s,n-s,K-s}$ is the output. That is, $Z$ is
the approximate solution to the relaxed problem (\ref{problem:relaxed}). Since $Z$ may not be binary valued,
our last task is to find a nearest matrix that is feasible in the original problem, i.e.~we
solve $\min_{X \in \vP_{m-s,n-s,K-s}} \|X-Z \|_F^2$.
Noting that $\|X-Z \|_F^2=\|X\|_F^2 -2 \TT Z^TX +\|Z\|_F^2$, we see that this is equivalent to
solving $\max_{X \in \vP_{m-s,n-s,K-s}} \TT Z^TX$, which is again a generalized linear assignment problem,
solvable efficiently, in time $O(\max\{m,n\}^3)$, via the Hungarian Algorithm, thanks to Lemma \ref{lemma:assign} (since the matrix 
in Hungarian Algorithm is $O(\max\{m,n\}$-by-$O(\max\{m,n\}$).
The optimal solution will be an approximate solution to the seeded subgraph-subgraph matching
problem as expressed in (\ref{sssmp}) and is the output of this algorithm which we call ssSGM. We summarize the
algorithm as follows:\\

\begin{algorithm}[H]
\linespread{1.45}\selectfont
\SetAlgoLined
\SetKwInOut{Input}{Input}
\SetKwInOut{Output}{Output}
\SetKwInOut{Initialization}{Initialization}
 \Input{$\pm 1$-Adjacency matrices $A \in \{ -1,1 \}^{m \times m}$, $B \in \{-1,1\}^{n \times n}$ of $G$,$H$ resp.,
 $K:=$ the number of core vertices, $s:=$ the number of seeds}
  \Initialization{$Z \leftarrow $ the matrix in $\vD_{m-s,n-s,K-s}$ with all entries $\frac{K-s}{(m-s)(n-s)}$ }
 \While{$Z$ is not approximately equal to what it was in immediately previous iteration}{
  Calculate $M \leftarrow A^{(21)}B^{(12)}+A^{(22)}ZB^{(22)} \in \R^{(m-s)\times (n-s)} $\;

  Compute $X^* \leftarrow \arg \max_{X \in \vP_{m-s,n-s,K-s}} \TT M^TX$ via Algorithm \ref{alg:glap} \;

  Solve $\max_{\lambda \in [0,1]} \TT A  ( I_s \oplus [\lambda X^*+(1-\lambda)Z]  ) B \left
  (I_s \oplus [\lambda X^*+(1-\lambda)Z] \right  )^T$; this is done exactly since the objective is
  quadratic in the one variable $\lambda$. Say $\lambda^*$ is the solution \;

  Set  $Z \leftarrow \lambda^* X^*+(1-\lambda^*)Z$ \;
  }

Compute $X^* \leftarrow \arg \max_{X \in \vP_{m-s,n-s,K-s}} \TT Z^TX$ via Algorithm \ref{alg:glap} \;

\Output{ sets $\Om := \{ v_i: \exists j \mbox{ s.t. } (I_s \oplus X^*)_{ij}=1 \},\OM := \{ w_j: \exists i \mbox{ s.t. } (I_s \oplus X^*)_{ij}=1 \}$. function $\phi :\Om \rightarrow \OM$ where  $\forall i,j$ \  $\phi (v_i)=w_j$ if and only if
         $(I_s \oplus X^*)_{ij}=1$ .}
\caption{ssSGM Algorithm  \label{alg:sssgm}}
\end{algorithm}
\vspace{.15in}

The ssSGM Algorithm (Algorithm \ref{alg:sssgm}) utilizes the following subroutine Algorithm \ref{alg:two} to solve
the generalized linear assignment problem; Algorithm \ref{alg:two} is correct because of 
Lemma \ref{lemma:assign}.
\vspace{.15in}

\begin{algorithm}[H]   \label{alg:two}
\linespread{1.45}\selectfont
\SetAlgoLined
\SetKwInOut{Input}{Input}
\SetKwInOut{Output}{Output}
\SetKwInOut{Initialization}{Initialization}
 \Input{Positive integers $c,d,e$ such that $e\leq \min \{c,d\}$ and matrix $M \in \R^{c \times d}$}

Calculate $\lb \leftarrow  \min \{ M_{ij} \}-1$, \ \ \ $\ub \leftarrow  \max \{ M_{ij} \}+1$\;

Construct $\mathcal{M}   \leftarrow \bigl [\begin{smallmatrix} M & \ub E \\ \ub E & \lb E \end{smallmatrix} \bigr ]
  \in \R^{(c+d-e)\times (c+d-e)} $,
  where matrix $E$ is all $1$s, appropriately sized   \;

Solve $\max_{\XX \in \vP_{c+d-e,c+d-e,c+d-e}}\TT \mathcal{M}^T \XX$ with Hungarian Algorithm \;

Define
  $X^*$  to be the first $c$ rows and first $d$ columns of the optimal $\XX$ in line 3   \;

\Output{$X^*\in \vP_{c,d,e}$ }
\caption{ Solving generalized
linear assignment problem:
$\max_{X \in \vP_{c,d,e}} \TT M^TX$    \label{alg:glap}}
\end{algorithm}
\vspace{.1in}

In practice, we restrict the number of Frank-Wolfe iterations (lines $1$ to $5$ in Algorithm \ref{alg:sssgm}) to $30$~iterations, thus the order of the runtime of ssSGM Algorithm (Algorithm \ref{alg:sssgm}) is the same order as the runtime of the Hungarian Algorithm in the subroutine Algorithm \ref{alg:glap}, which is  $O(\max\{m,n\}^3)$.

In the special case where $m=n=K$, the ssSGM algorithm is the SGM algorithm of \citet{sgm}, the one difference being that SGM can and does use the standard
$\{ 0,1\}$-adjacency matrices instead of the $\pm 1$-adjacency matrices used in ssSGM; indeed, ssSGM can't use $\{ 0,1\}$-adjacency matrices in general.

\section{Experiments with synthetic data}
\label{sec:exp}

In this section we perform numerical experiments to illustrate the utility of the ssSGM algorithm.

We will be working with the correlated Erdos-Renyi random graph model, a special case
of the correlated Bernoulli random graph model in Section \ref{section:model}.
The model parameters are positive integers $m,n,K$ such that $K \leq \min \{m,n\}$ and real numbers $p,\varrho \in [0,1]$.
The vertex set of $G$ is $V=\{v_1,v_2,\ldots,v_{m}\}$ and the vertex set of $H$ is $W=\{w_1,w_2,\ldots, w_{n}\}$.
For all positive integers $i,j$ such that $1\leq i < j \leq m$, \ $p$ is the probability that
$v_i \sim_{G} v_j$ and, for all positive integers 
$i,j$ such that $1 \leq i< j \leq n$,
$p$ is the probability that $w_i \sim_{H} w_j$.
For all positive integers $i,j$ such that $1\leq i < j \leq K$, \  $\varrho$ is the Pearson correlation coefficient for the Bernoulli random variable 
indicators ${\bf 1}_{v_i \sim_G v_j}$ and 
${\bf 1}_{w_i \sim_H w_j}$;
all other adjacencies are collectively independent. The random graphs $(G,H)$ are
called {\it correlated Erdos-Renyi random graphs}. The parameters $m,n,K,p,\varrho$  fully specify the distribution of $(G,H)$.
Clearly, the core here are the vertices $\{v_1,v_2,\ldots,v_K\}$ and  $\{w_1,w_2,\ldots,w_K\}$, with the natural alignment
associating $v_i$ to $w_i$, for all $i\in [K]$.
There are $s$ seeds, for a given nonnegative integer $s \leq K$, and they are $\{v_1,v_2,\ldots,v_s\}$ and $\{w_1,w_2,\ldots,w_s\}$.

After realizing $G$ and $H$ from this model, we discrete-uniformly randomly relabel
the
nonseed vertices $\{v_{s+1},v_{s+2},\ldots,v_{m} \}$ of $G$
and the
nonseed vertices $\{w_{s+1},w_{s+2},\ldots,w_{n} \}$ of $H$.
We use the ssSGM algorithm with the $s$ seeds to attempt to identify the $K$ core vertices in $G$ and
in $H$ and to align them; the ssSGM's {\it match ratio}
is the fraction of the $K-s$  nonseed, core vertices of $G$ that the algorithm correctly matches to their naturally aligned counterpart
in $H$. 

If $m=n$ then we also use the SGM algorithm to do a seeded graph match of $G$ to $H$ (of course, if $m\ne n$ then SGM does not work).
Recall that the SGM algorithm applies Frank-Wolfe methodology to match the entire graphs, and not just
subgraphs, which is the innovation of ssSGM. As above, SGM's {\it match ratio} is the fraction
of the $K-s$  nonseed, core vertices of $G$ that the algorithm correctly matches.

When running the SGM algorithm, $\pm 1$-adjacency matrices were used instead of the usual $\{ 0,1\}$-adjacency matrices. This
was done in order to get the most fair comparison with the ssSGM algorithm, which uses $\pm 1$-adjacency matrices and doesn't work 
with $\{ 0,1\}$-adjacency matrices. However, empirical experiments indicate that the choice of $\pm 1$-adjacency matrices v.s.~$\{ 0,1\}$-adjacency matrices in the SGM
algorithm doesn't significantly impact the effectiveness of the SGM algorithm. See \ref{appendixD} for an empirical demonstration of this.

\begin{figure}[hp!]
     \centering
     \begin{subfigure}[b]{0.32\textwidth}
         \centering
         \includegraphics[width=\textwidth]{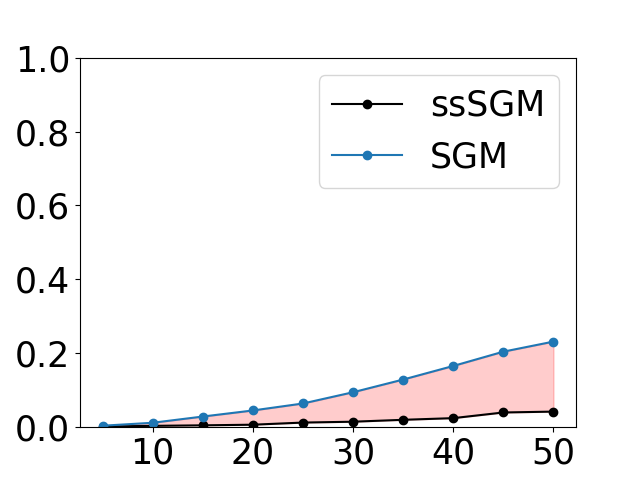}
         \caption*{$ p = 0.1, \varrho = 0.5$}
     \end{subfigure}
     \hfill
     \begin{subfigure}[b]{0.32\textwidth}
         \centering
         \includegraphics[width=\textwidth]{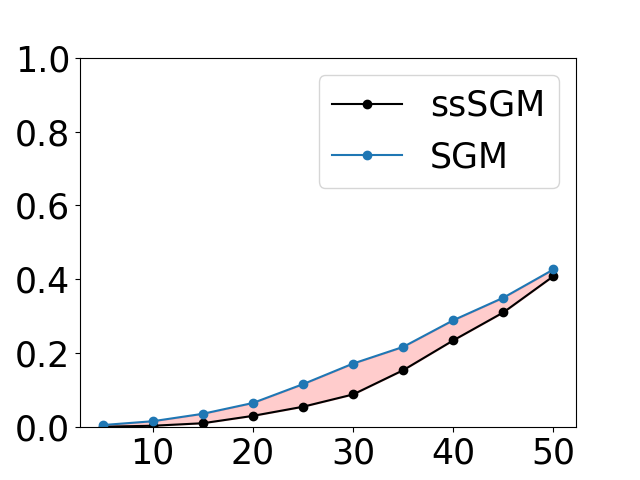}
         \caption*{$p = 0.3,\varrho = 0.5 $}
     \end{subfigure}
     \hfill
     \begin{subfigure}[b]{0.32\textwidth}
         \centering
         \includegraphics[width=\textwidth]{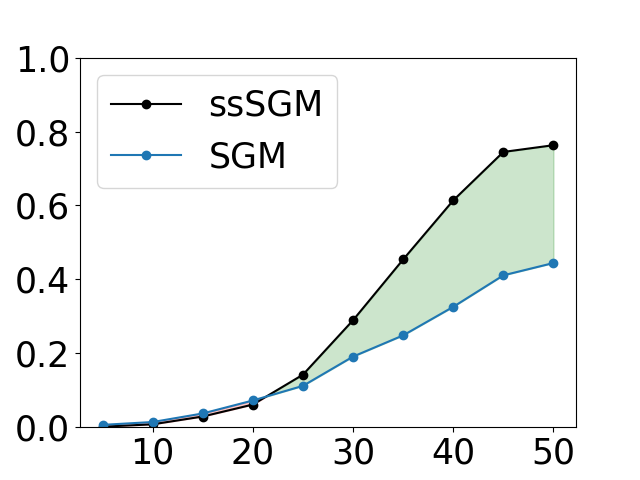}
         \caption*{$ p = 0.5, \varrho = 0.5$}
     \end{subfigure}
     \\
     \centering
     \begin{subfigure}[b]{0.32\textwidth}
         \centering
         \includegraphics[width=\textwidth]{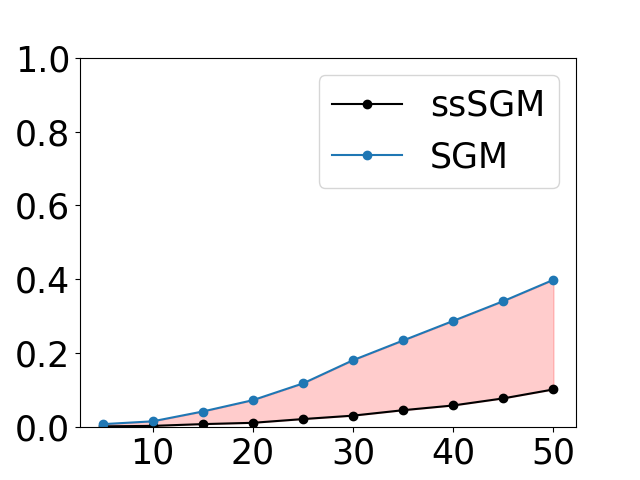}
         \caption*{$ p = 0.1, \varrho = 0.6$}
     \end{subfigure}
     \hfill
     \begin{subfigure}[b]{0.32\textwidth}
         \centering
         \includegraphics[width=\textwidth]{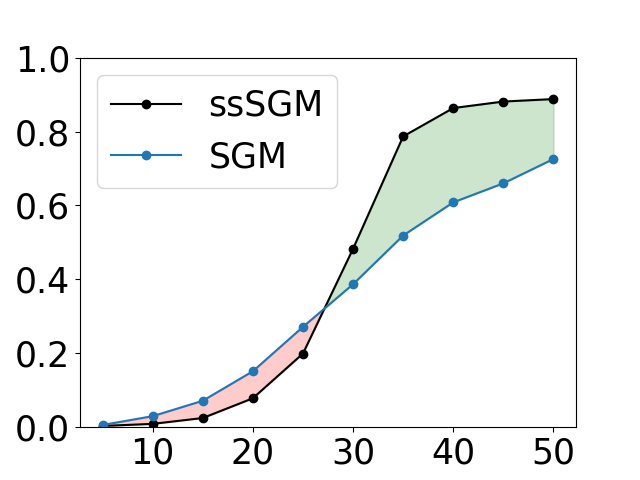}
         \caption*{$ p = 0.3, \varrho = 0.6$}
     \end{subfigure}
     \hfill
     \begin{subfigure}[b]{0.32\textwidth}
         \centering
         \includegraphics[width=\textwidth]{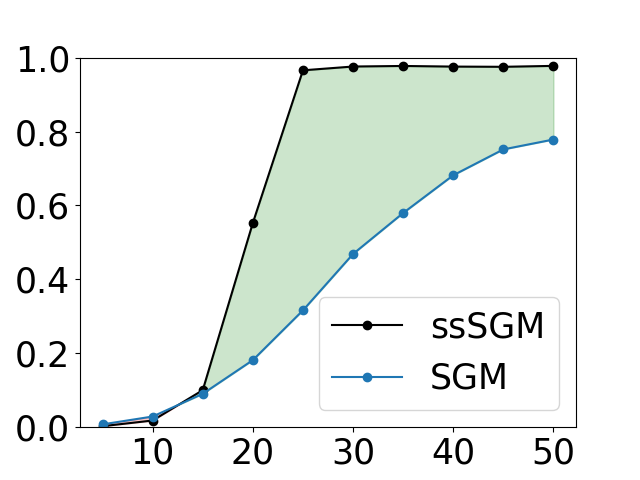}
         \caption*{$ p = 0.5, \varrho = 0.6$}
     \end{subfigure}
     \\
     \centering
     \begin{subfigure}[b]{0.32\textwidth}
         \centering
         \includegraphics[width=\textwidth]{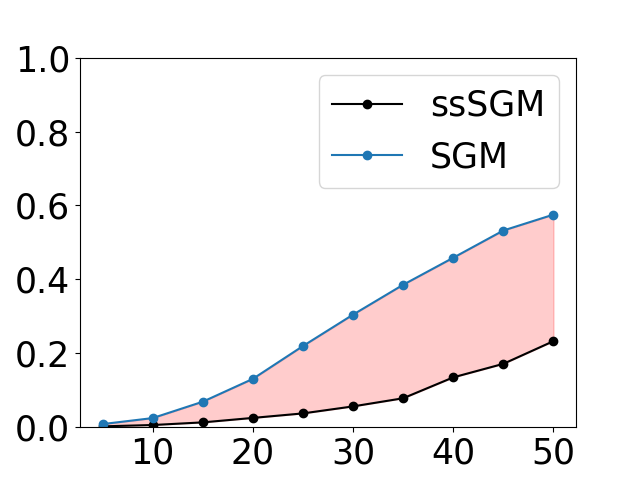}
         \caption*{$ p = 0.1, \varrho = 0.7$}
     \end{subfigure}
     \hfill
     \begin{subfigure}[b]{0.32\textwidth}
         \centering
         \includegraphics[width=\textwidth]{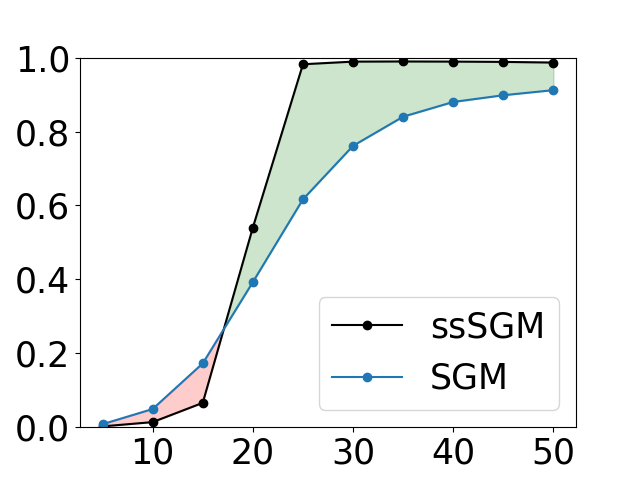}
         \caption*{$ p = 0.3, \varrho = 0.7$}
     \end{subfigure}
     \hfill
     \begin{subfigure}[b]{0.32\textwidth}
         \centering
         \includegraphics[width=\textwidth]{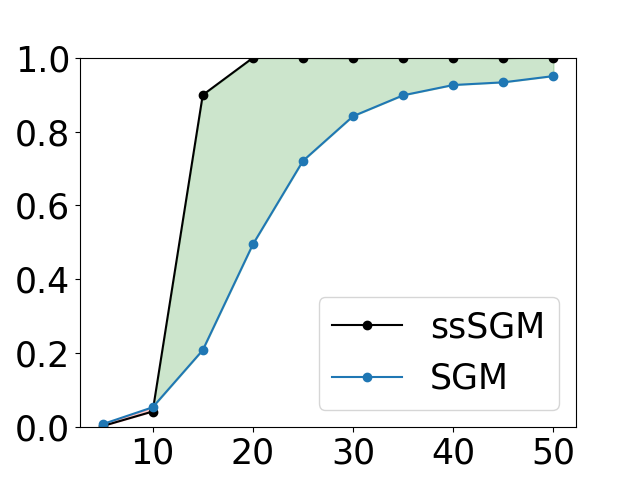}
         \caption*{$ p = 0.5, \varrho = 0.7$}
     \end{subfigure}
     \\
     \centering
     \begin{subfigure}[b]{0.32\textwidth}
         \centering
         \includegraphics[width=\textwidth]{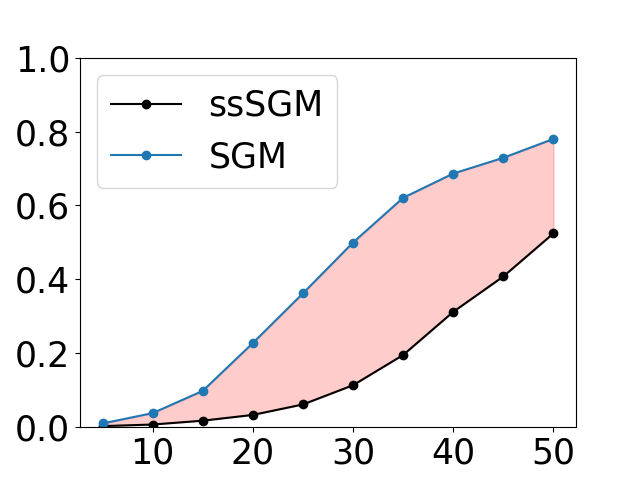}
         \caption*{$ p = 0.1, \varrho = 0.8$}
     \end{subfigure}
     \hfill
     \begin{subfigure}[b]{0.32\textwidth}
         \centering
         \includegraphics[width=\textwidth]{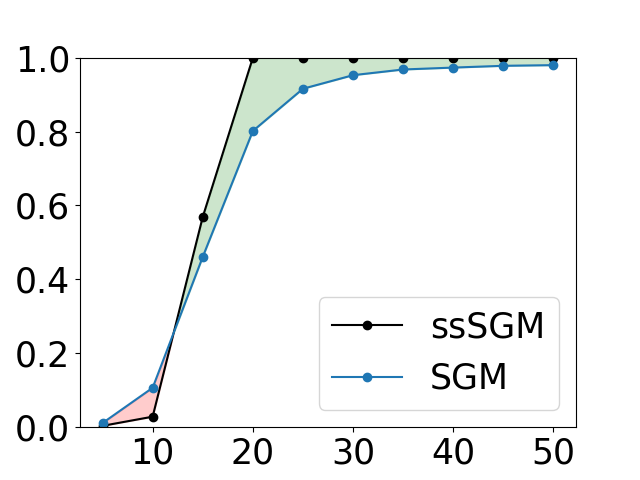}
         \caption*{$ p = 0.3, \varrho = 0.8$}
     \end{subfigure}
     \hfill
     \begin{subfigure}[b]{0.32\textwidth}
         \centering
         \includegraphics[width=\textwidth]{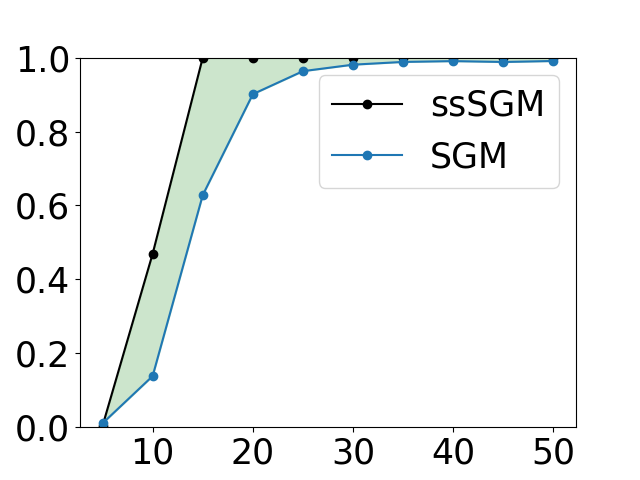}
         \caption*{$p = 0.5, \varrho = 0.8$}
     \end{subfigure}
    \caption{Comparing match ratio of ssSGM vs SGM on synthetic data, $m=n=500$, $K=100$ for different combinations of the values of $p$ and $\varrho$.
    {\bf In each of these figures, the ``$x$" axis is labeled with the number of seeds, and the ``$y$" axis is labeled with the match ratio.}\label{fig:synth}}
\end{figure}

The experiments here had $m=n=500$, $K=100$, and
we did $200$ experiments for each combination of $p=0.1, 0.3 ,0.5$ and $\varrho=0.5, 0.6, 0.7, 0.8$ and
$s=5,10,15,20,\ldots, 50$; for each such combination, we averaged the match ratio for each of ssSGM and SGM over the $200$ experiments.
Figure \ref{fig:synth} plots the average match ratios of ssSGM and SGM against the number of seeds $s$, for each combination of $p$ and $\varrho$. The
difference between the average match ratios of ssSGM and SGM is highlighted in green when the average match ratio of ssSGM was higher than that of SGM, and
the difference is highlighted in red otherwise.

Note that ssSGM is empirically seen here to achieve a better match ratio than SGM as the graphs gets less sparse (less dense, if the density is more than $0.5$) and as the correlation between the two graphs gets stronger. However, ssSGM identifies the core vertices in the two graphs, which SGM is not designed to do.
So, even in cases where SGM achieves a better match ratio on the core, this might be useless as a practical matter if the core vertices need to be identified. Indeed, {\it it is important to emphasize the ssSGM has a more substantial inferential task than SGM}.

Also note that, in these experiments, higher correlation and greater density led to a better match ratio in general for ssSGM (not just in comparing ssSGM to SGM). Indeed, in Section \ref{section:theo}, we will state and prove a consistency result regarding 
matchability for subgraph-subgraph matching in which, to guarantee matchability, (say $0<p<\frac{1}{2}$, $m\leq n$, and $s=0$) needs $\varrho > (1-2p)^2$ and $\left ( \varrho - (1-2p)^2 \right )^2 K \in \omega (\log n)$, see Section \ref{section:theo} for critical details. These experiments reflect those provisions; notably, guaranteed matchability is possible if $p$ is close to $0$ (of order $\omega(\sqrt{(\log n)/K})$) as long as the correlation is sufficiently close to $1$. 

\begin{figure}[hp!]
    \centering
    \includegraphics[width=0.32\linewidth]{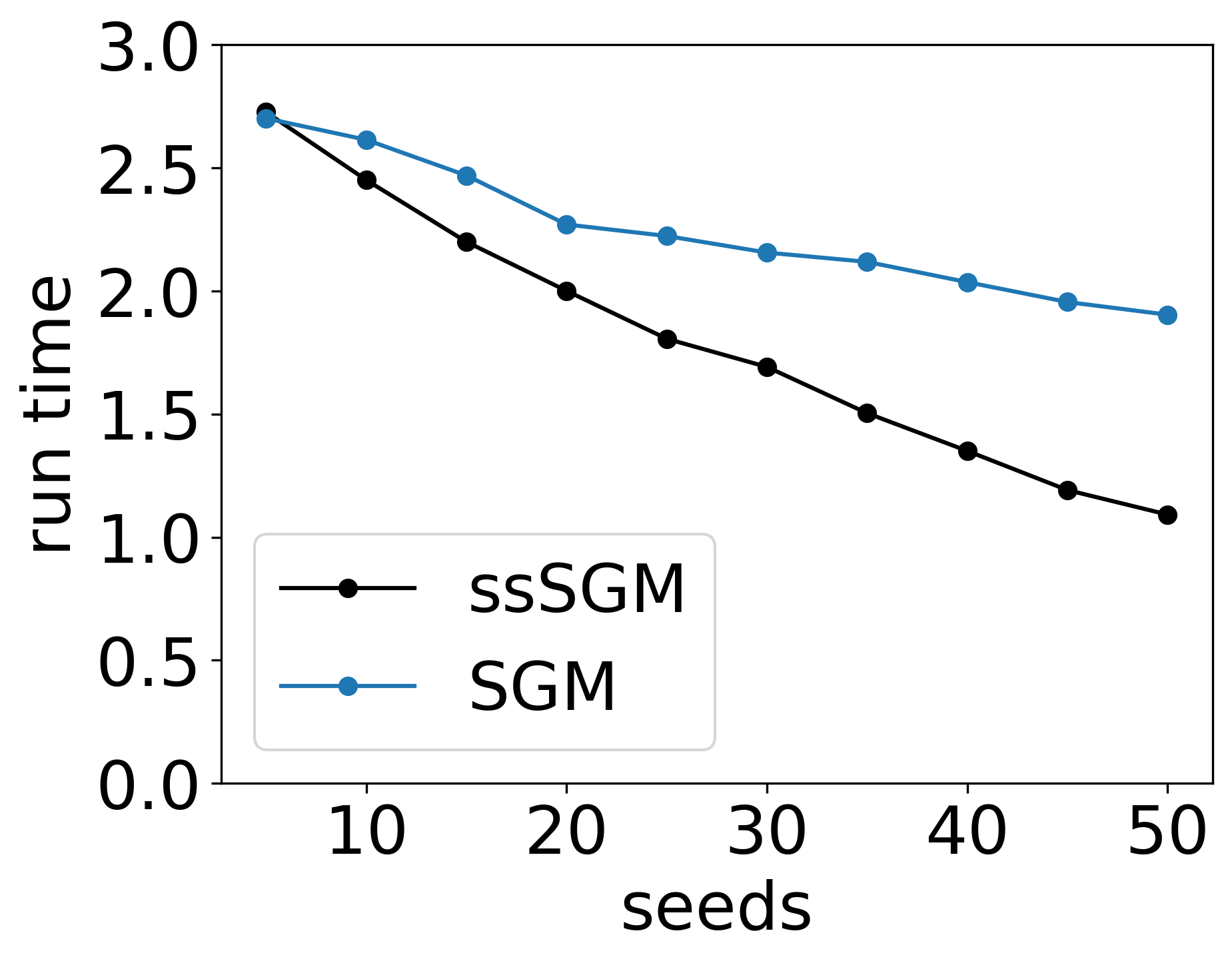}
    \includegraphics[width=0.32\linewidth]{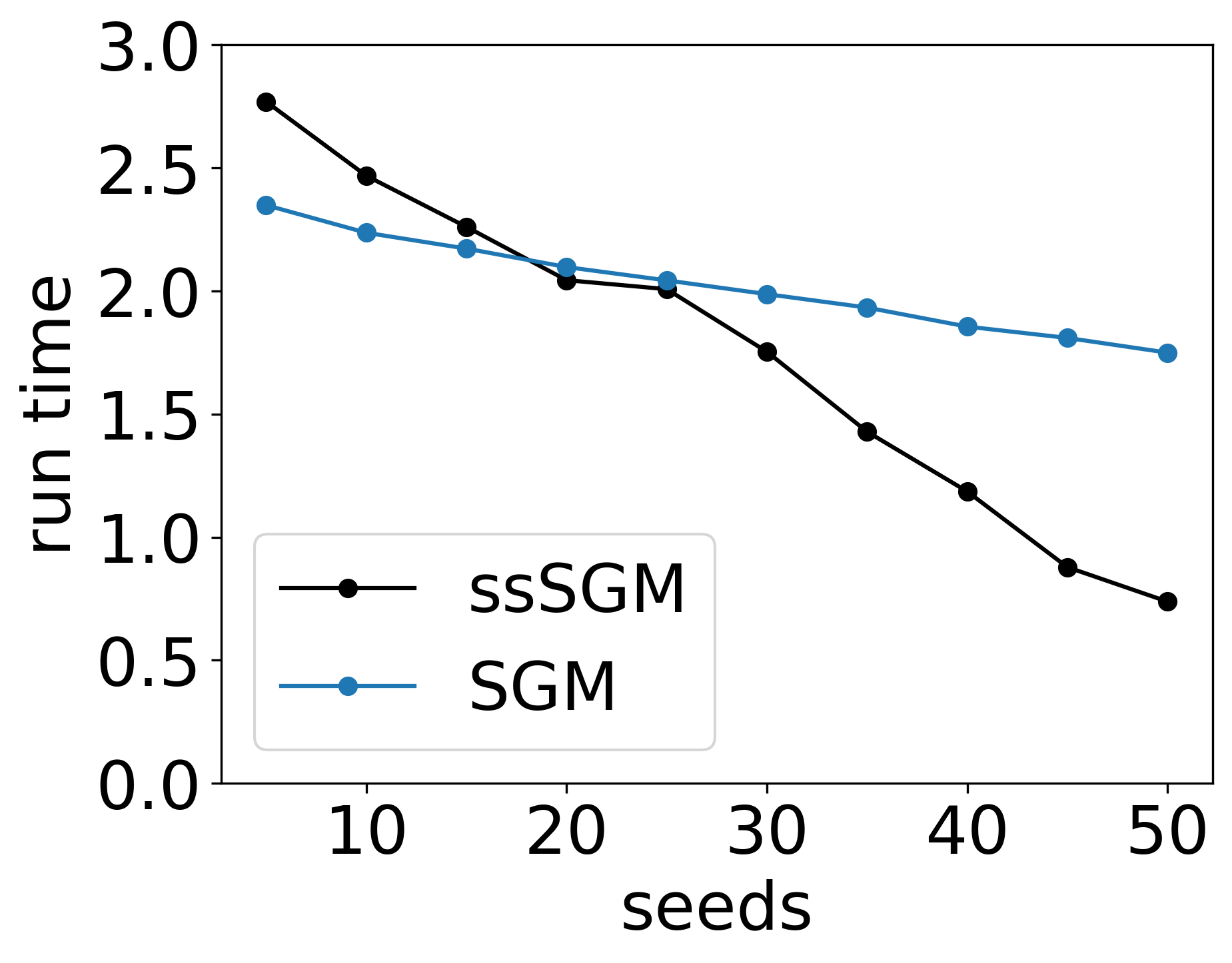}
    \includegraphics[width=0.32\linewidth]{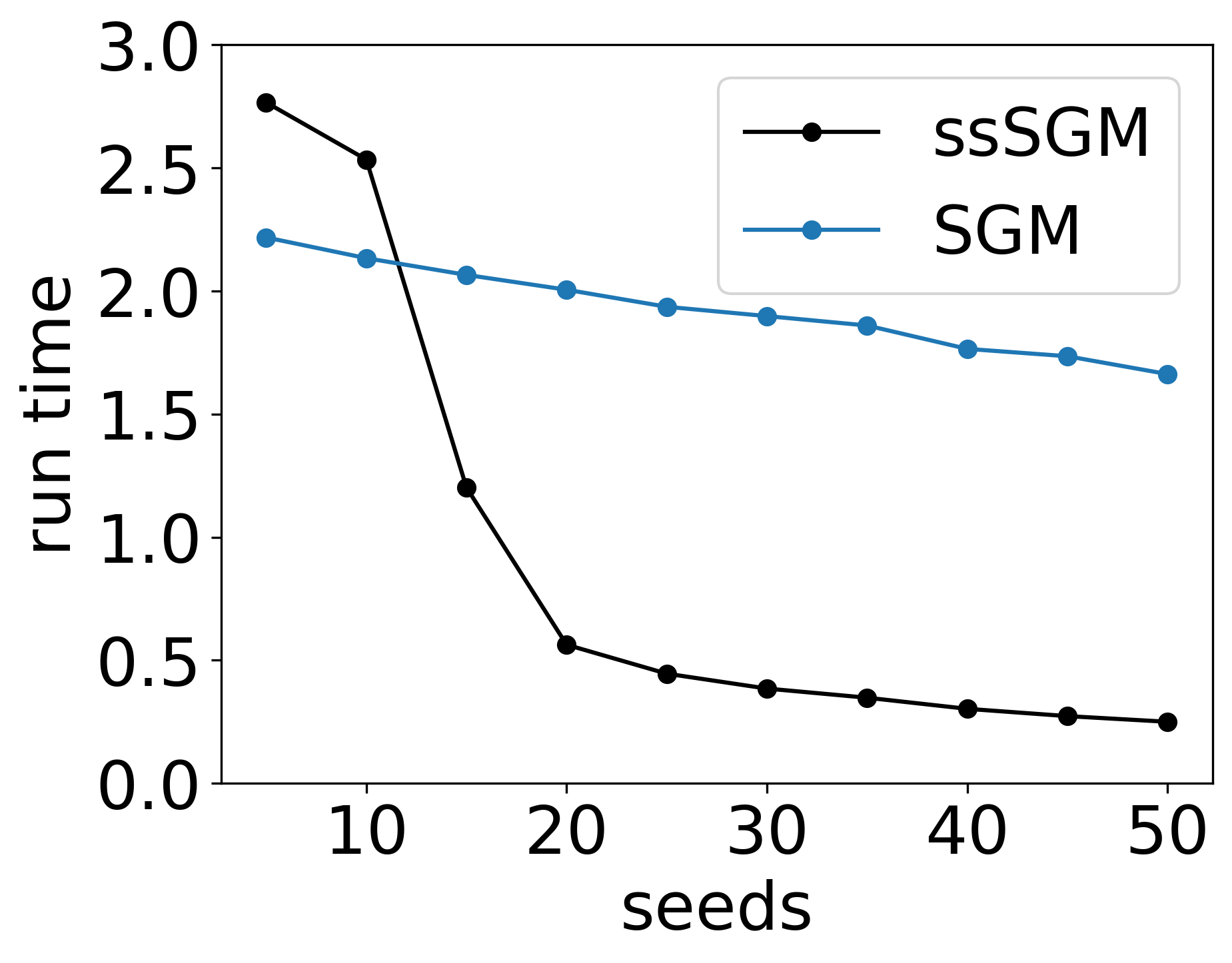}
    \caption{Average runtime in seconds of ssSGM and SGM as a function of the number of seeds $s$ for $p = 0.1, \varrho = 0.5$ (left panel), $p = 0.3, \varrho = 0.6$ (middle panel), $p = 0.5, \varrho = 0.7$ (right panel)}
    \label{fig:runtime}
\end{figure}

In Figure \ref{fig:runtime} we plot the average runtime in seconds (on a MacBook Pro, M2, 2022 running code in Python) of ssSGM and SGM as a function of the number of seeds $s$, but just for the above experiments where $p=0.1,\varrho=0.5$ (left panel),
then where $p=0.3, \varrho=0.6$ (middle panel),
then where $p=0.5,\varrho=0.7$ (right panel).

\begin{figure}
    \centering
    \begin{subfigure}[b]{0.47\textwidth}
    \centering
    \includegraphics[width=\textwidth]{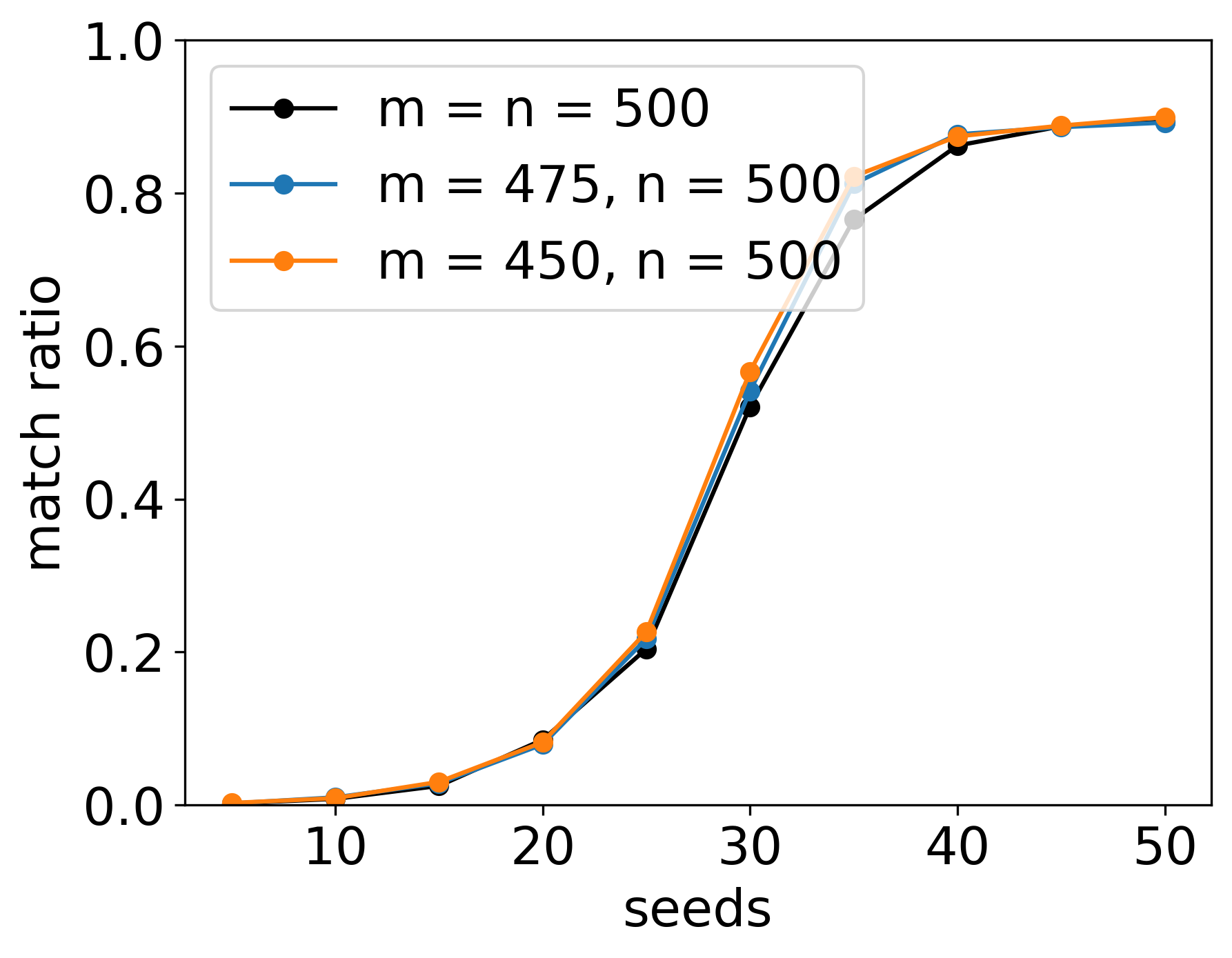}
    \caption*{$ p = 0.3, \varrho = 0.6$}
    \end{subfigure}
    \hfill
    \begin{subfigure}[b]{0.47\textwidth}
    \centering
    \includegraphics[width=\textwidth]{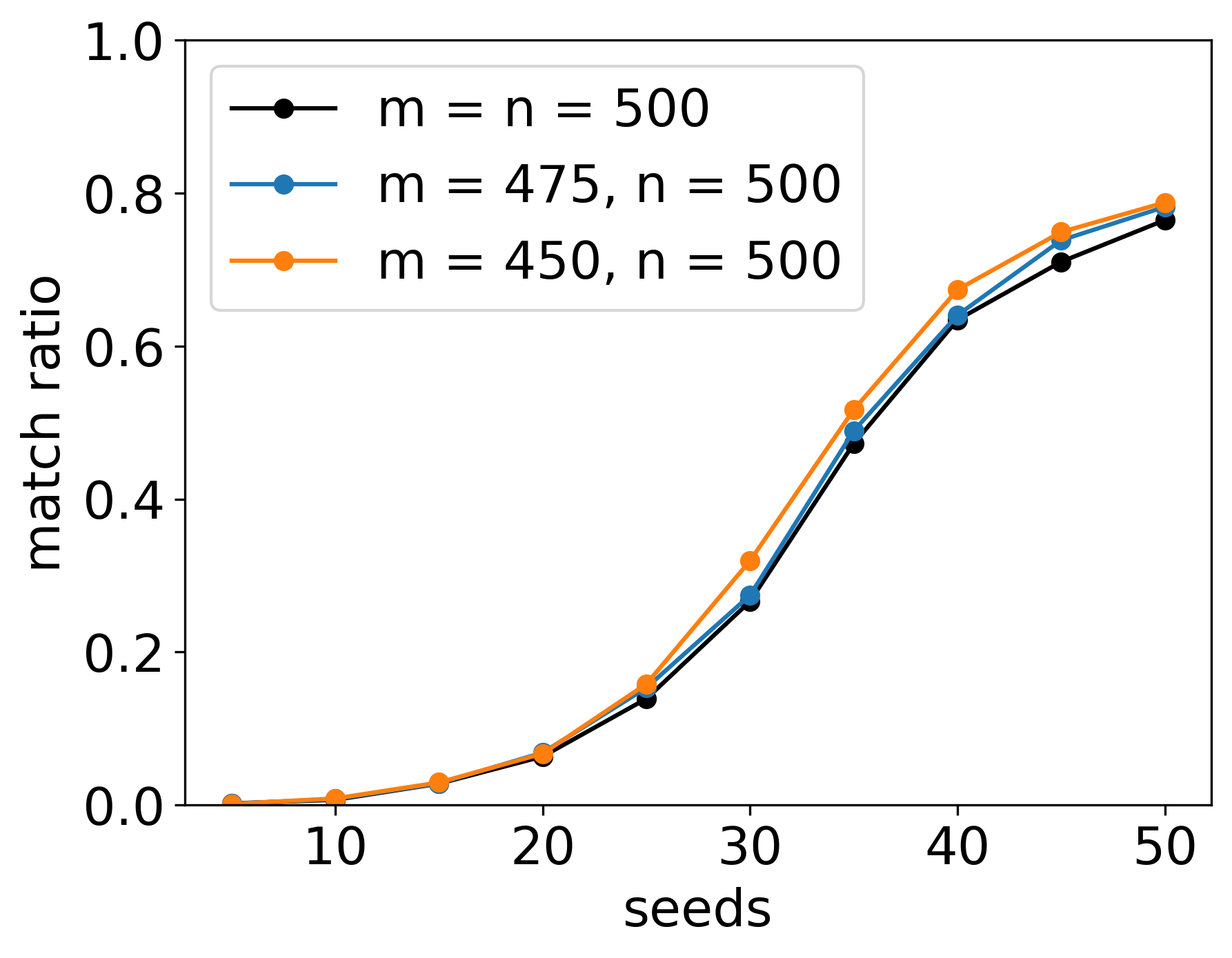}
    \caption*{$ p = 0.5, \varrho = 0.5$}
    \end{subfigure}
    \caption{Match ratio as a function of number of seeds for ssSGM, comparing $m=n$ to $m\neq n$.}
    \label{fig:diff-size}
\end{figure}

The above experiments were done with $m=n$ in order to compare ssSGM with SGM. For $p=0.3,\varrho=0.6$ and again for $p=0.5, \varrho=0.5$, we repeated the above experiments, but  this time with $m=450$, $n=500$, $K=100$ and then again with 
and $m=475, n=500$, $K=100$. 
Figure \ref{fig:diff-size}  plots the match ratio of ssSGM as a function of the number of seeds $s$; the left panel for 
$p=0.3$, $\varrho=0.6$ for each of $m=450,475,500$ (when $m=500$ this was the previous experiments), and the right panel for 
$p=0.5$, $\varrho=0.5$ for each of $m=450,475,500$ (when $m=500$ this was the previous experiments). It is seen here that ssSGM was not impeded by the different sized graphs.

In \ref{appendixC} we perform additional experiments in the context of stochastic block models and correlated Bernoulli random graph models to demonstrate the effectiveness of ssSGM.

\section{Experiments with real data}
\label{sec:real}

The first data set for this section comes from the University of California Irvine Machine Learning Repository.
``Mathwiki" is a simple graph with $1068$ vertices and $27079$ edges; the vertices are Wikipedia pages for mathematics topics and, for
every pair of vertices, the two vertices are adjacent precisely when one associated webpage has a hyperlink to the other associated webpage. The data set is found at \url{https://archive.ics.uci.edu/ml/datasets/Wikipedia+Math+Essentials}.

For each value of $K=700, 800, 900, 1000$ and each number of seeds $s = 20,40,60,...,500$, we perform the following experiment $200$ times.
We randomly select $K$ core vertices from the $1068$ vertices of Mathwiki. The set $V$ consists of the $K$ core vertices and $\frac{1068-K}{2}$ randomly selected noncore vertices of Mathwiki, and the set $W$ consists of the $K$ core vertices and the
other $\frac{1068-K}{2}$ noncore vertices of Mathwiki.
The graph $G$ is taken to be the subgraph of Mathwiki induced by $V$, and the graph $H$ is taken to be the subgraph of Mathwiki induced by $W$. The $K$ core vertices are common to both $G$ and $H$ with the identity mapping, and we randomly choose $s$ seeds from the $K$ core vertices. The nonseed vertices of $G$ are randomly permuted and the nonseed vertices of $H$ are randomly permuted.
We run ssSGM and SGM on $(G,H)$, and average the core match ratios for ssSGM and for SGM over the $200$ repetitions.

  \begin{figure}[h!]
     \centering
     \begin{subfigure}[b]{0.48\textwidth}
         \centering
         \includegraphics[width=\textwidth]{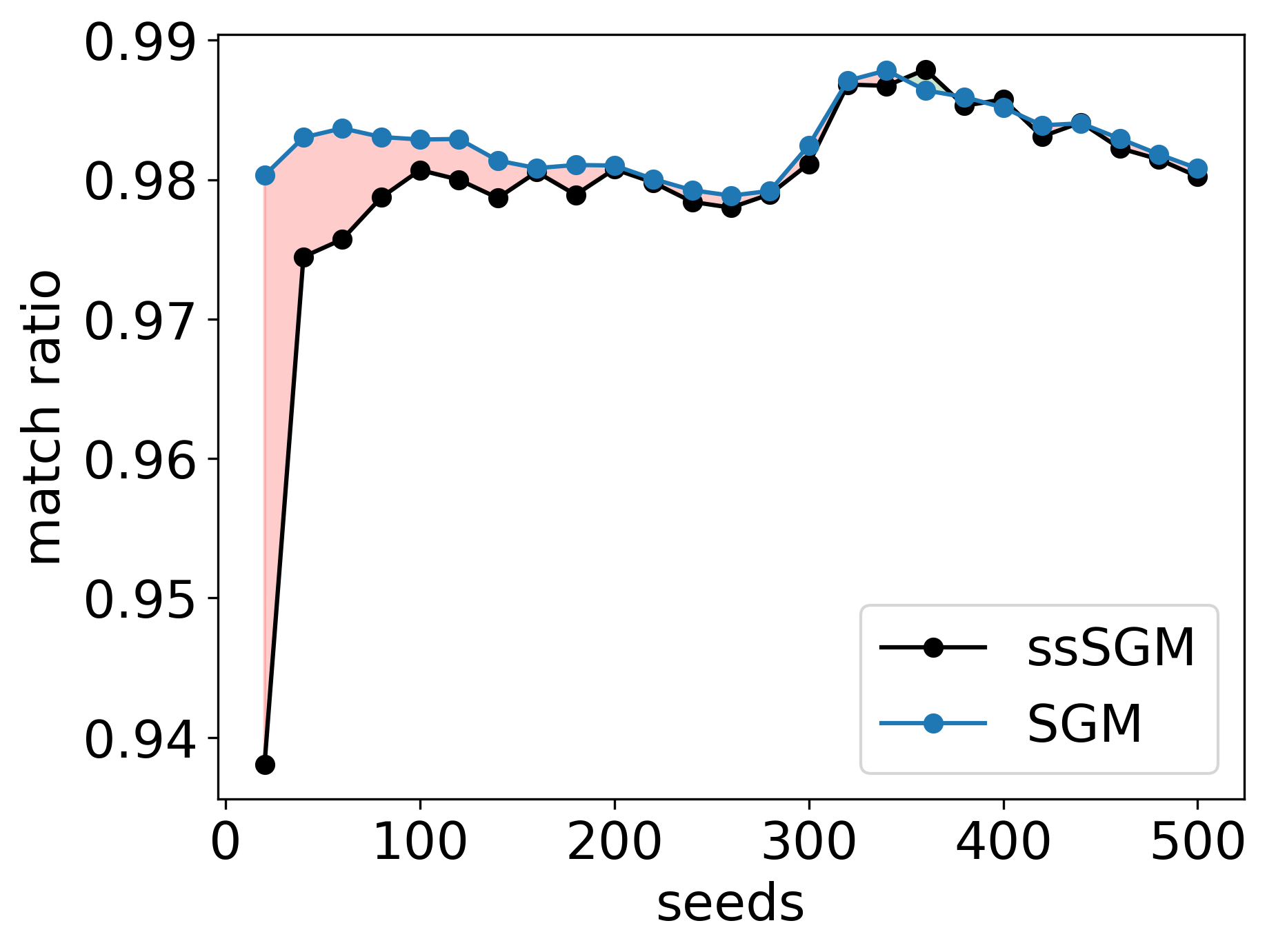}
         \caption*{$K = 700,\ \  m = n = 884$}
     \end{subfigure}
     \hfill
     \begin{subfigure}[b]{0.48\textwidth}
         \centering
         \includegraphics[width=\textwidth]{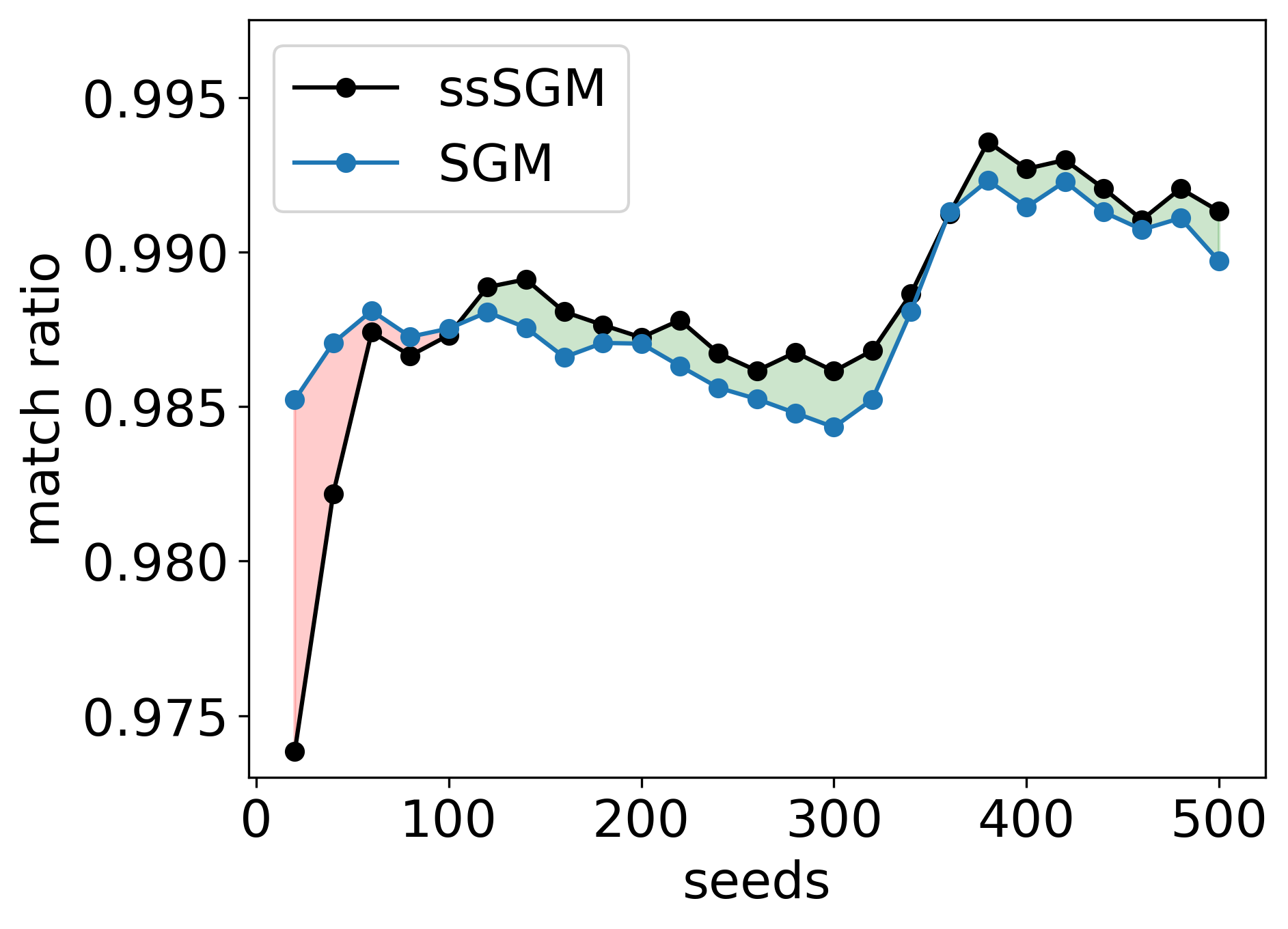}
         \caption*{$K=800, \ \  m = n = 934$}
     \end{subfigure}

     \centering
     \begin{subfigure}[b]{0.48\textwidth}
         \centering
         \includegraphics[width=\textwidth]{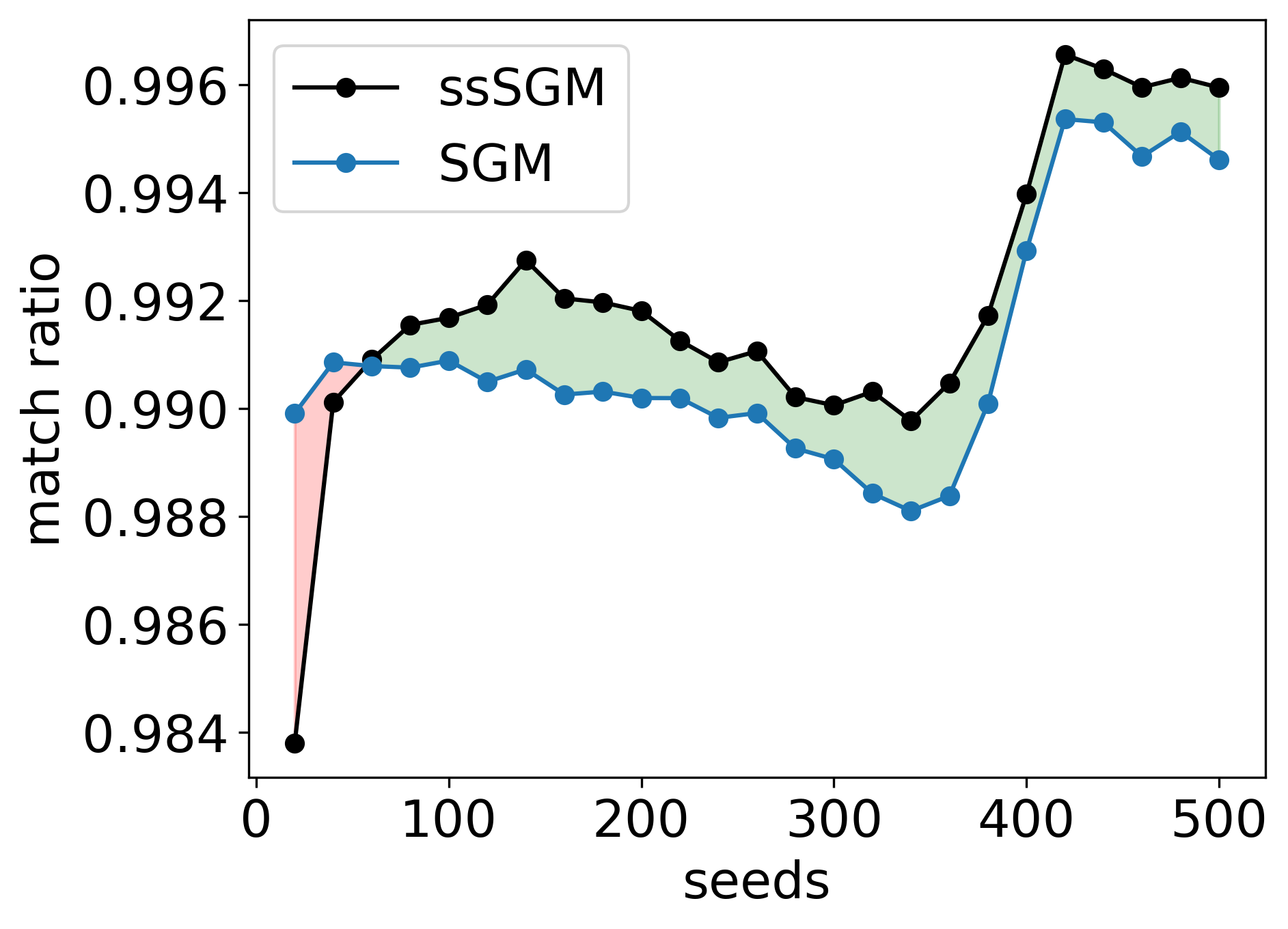}
         \caption*{$K = 900, \ \ m = n = 984$}
     \end{subfigure}
     \hfill
     \begin{subfigure}[b]{0.48\textwidth}
         \centering
         \includegraphics[width=\textwidth]{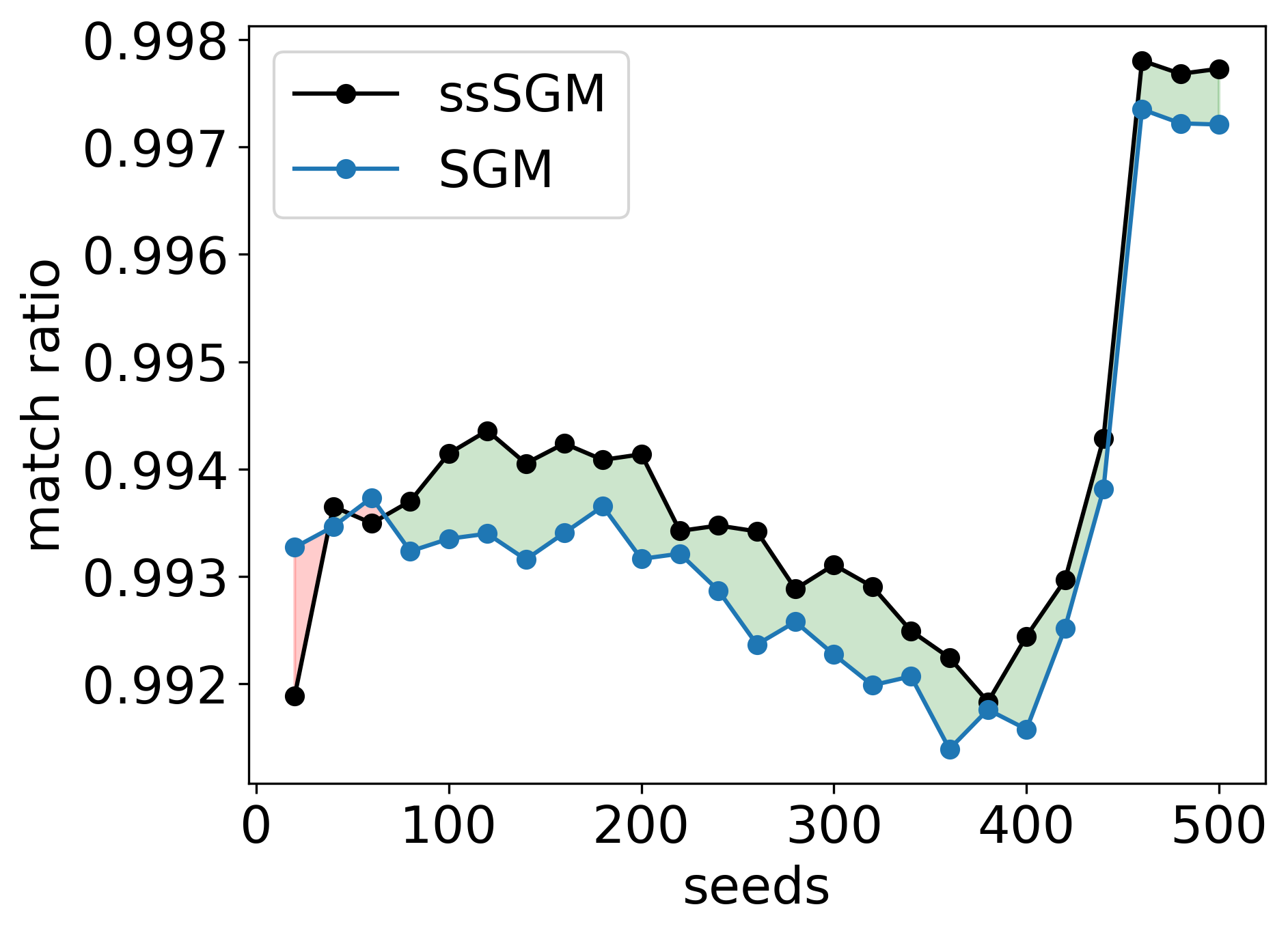}
         \caption*{$K = 1000,  \ \ m = n = 1034$}
     \end{subfigure}
     \caption{Mathwiki experiments \label{fig:real}}
\end{figure}

Figure \ref{fig:real} plots the average match ratios of ssSGM and SGM against the number of seeds $s$, separately for each value of $K$.
The
difference between the average match ratios of ssSGM and SGM is highlighted in green when the average match ratio of ssSGM was higher than that of SGM, and
the difference is highlighted in red otherwise.

Real data is often quite sparse, and sometimes the correlation between core vertices can be quite weak. In those circumstances, SGM can achieve a higher core match ratio than ssSGM but, as mentioned before, SGM would be useless if the core needs to be identified, a task that ssSGM is designed to perform and SGM is not. The Mathwiki graph has density $\frac{27079}{\binom{1068}{2}}=0.0475$, which is somewhat sparse, but the
core correlation is perfect; indeed, the match ratios were high, and ssSGM often outperformed SGM in average match ratio. Note that the $y$-axis in
Figure \ref{fig:synth} was from $0$ to $1$, whereas the $y$-axis in Figure \ref{fig:real} is scaled much closer to $1$.

We next perform another set of experiments with real data that we collected, similar to the above experiments, but the two graphs have a different number of vertices. 
The graph $G$ has vertex set $V$ consisting of all the Wikipedia pages hyperlinked from the 
Wikipedia page ``{\textbf{Mathematics}."
The graph $H^{(1)}$ has vertex set $W^{(1)}$ consisting of all the Wikipedia pages hyperlinked from the Wikipedia page      ``{\textbf{Geometry}."
The graph $H^{(2)}$ has vertex set $W^{(2)}$ consisting of all the Wikipedia pages hyperlinked from the Wikipedia page      ``{\textbf{Algebra}." For each 
of these three graphs, the edge 
sets are the pairs of vertices hyperlinked from either one to the other. These graphs were obtained on February 20, 2025; 
$G$ has  $746$ vertices and $23950$ edges, $H^{(1)}$ has $605$ vertices and $23160$ edges, and $H^{(2)}$ has $464$ vertices and $15905$ edges.
For $G$ and $H^{(1)}$, the core 
$V \cap W^{(1)}$ has $K=|V \cap W^{(1)}|=201$ vertices. 
For $G$ and $H^{(2)}$, the core 
$V \cap W^{(2)}$ has $K=|V \cap W^{(2)}|=161$ vertices.
Because the number of vertices in these graphs are not the same we are unable to~run~SGM.

For each number of seeds $s=10,20,30,\ldots,160$, we did $200$ experiments of randomly sampling $s$ seeds from the core of $G$ and $H^{(1)}$ and performed ssSGM on $G$ and $H^{(1)}$;
the left panel of Figure \ref{fig:real2} plots
the average match ratio vs the number of seeds $s$.
For each number of seeds $s=10,20,30,\ldots,120$, we did $200$ experiments of randomly sampling $s$ seeds from the core of $G$ and $H^{(2)}$ and performed ssSGM on $G$ and $H^{(2)}$;
the right panel of Figure \ref{fig:real2} plots
the average match ratio vs the number of seeds $s$.
\begin{figure}[h!]
\centering
\includegraphics[width=0.48\textwidth]{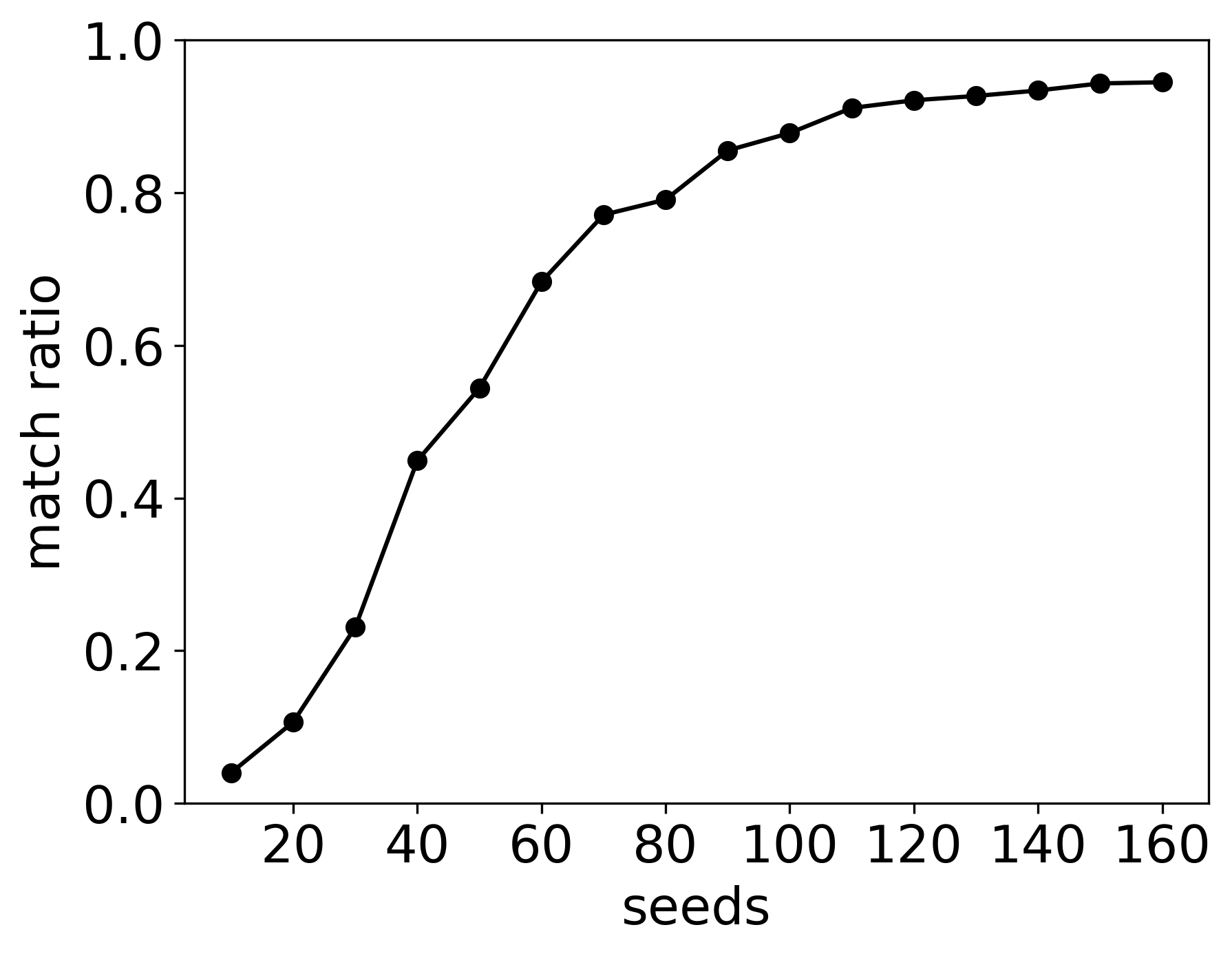}
\includegraphics[width=0.48\textwidth]{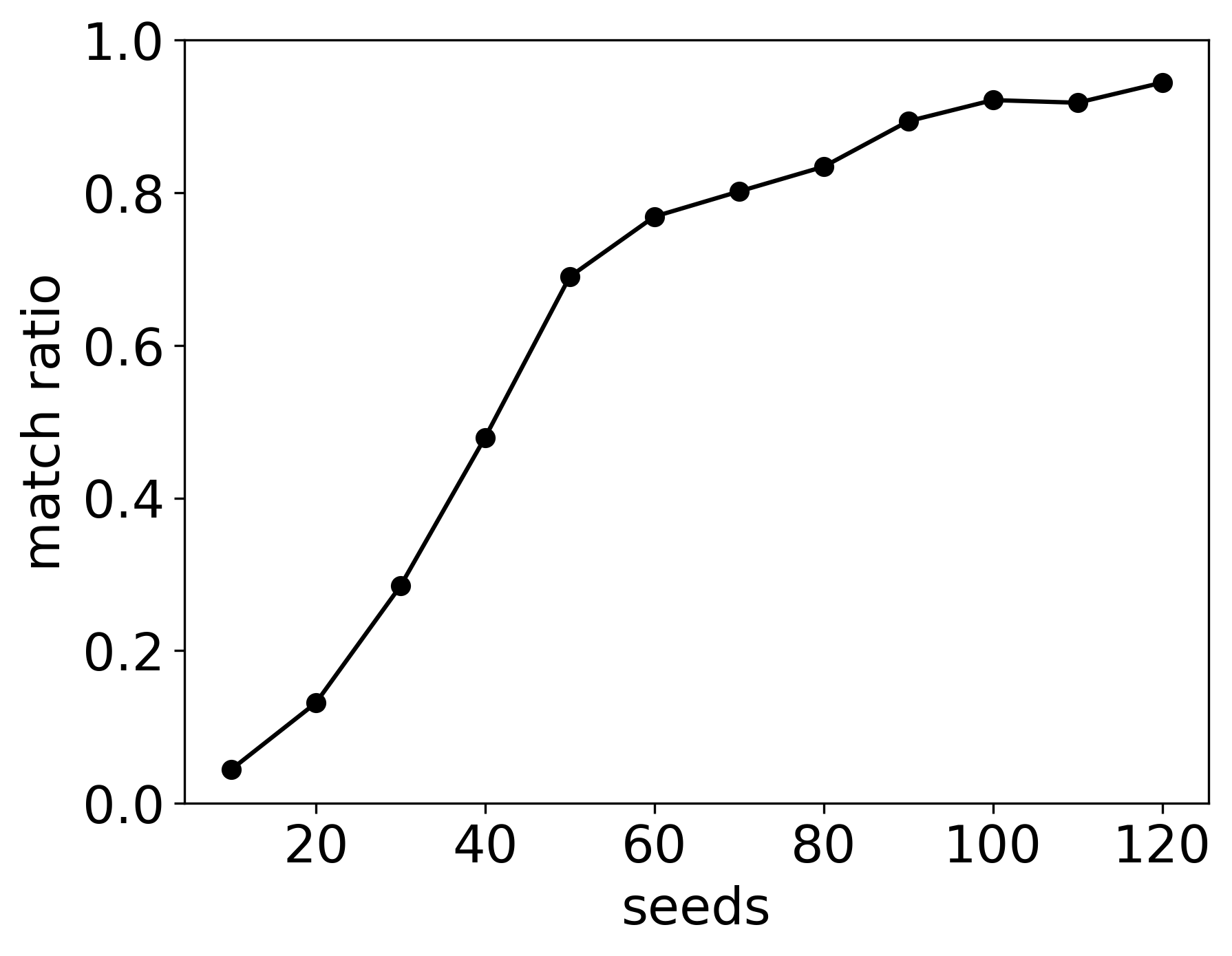}
\caption{Web crawl experiments for Mathematics vs Geometry (left) and for Mathematics vs Algebra (right).} \label{fig:real2}
\end{figure}

\section{Consistency \label{section:theo}   }

Until this point in this paper, our work has been focused on developing the ssSGM algorithm and highlighting it as an efficient and effective tool for an approximate solution to the seeded subgraph-subgraph matching problem. This section, Section \ref{section:theo},  is concerned with matchability from the information theory direction and is in the spirit of \citet{sussman2018matched,rel}. Specifically, when there is a natural underlying correspondence between a
core of a pair of graphs, when will the (exact, not just approximate) solution to the seeded subgraph-subgraph matching problem identify this core and give us this natural correspondence?
Our main result for this work is Theorem \ref{thm:main} which, under the model described in Section \ref{section:model}, guarantees that this almost always happens when there are mild assumptions in place.

\subsection{The random graph model and the main consistency result \label{section:model} }

We first define the correlated Bernoulli random graph model with which we will work. For any positive integer $k$, recall that
$[k]$ denotes the set of integers $\{1,2,3,\ldots,k\}$.
The correlated Bernoulli random graph model parameters are:
\begin{itemize} 
\item Positive integers $m,n,K$ such that $K \leq \min \{ m,n \}$,
\item Real numbers $p_{ij} \in [0,1]$, for all positive integers $i,j$ such that $1\leq i < j \leq \max \{ m,n \}$, and 
\item Real numbers $\varrho_{ij} \in [0,1]$, for all positive integers $i,j$ such that $1\leq i < j \leq K$. 
\end{itemize}
For simplicity,
the vertex set of $G$ is $V=[m]$ and the vertex set of $H$ is $W=[n]$.
For all positive integers $i,j$ such that $1\leq i < j \leq m$, \  $p_{ij}$ is the probability that $i \sim_{G} j$ and,
for all positive integers $i,j$ such that $1\leq i < j \leq n$, \ $p_{ij}$ is the probability that $i \sim_{H} j$
and, for all positive integers $i,j$ such that $1\leq i < j \leq K$, \  $\varrho_{ij}$ is the Pearson correlation coefficient for
the Bernoulli random variable 
indicators ${\bf 1}_{i \sim_G j}$ and 
${\bf 1}_{i \sim_H j}$;
all other adjacencies are collectively independent. The random graphs $(G,H)$ are
called {\it correlated Bernoulli random~graphs}.

The parameters fully specify the distribution of $(G,H)$.
Clearly, the core here are the vertices $[K] \subseteq V$ and  $[K] \subseteq W$, with the natural alignment
associating $i \in V$ to $i \in W$, for all $i\in [K]$.
  For our main result in this section, we will not need to suppose there are seeds.\\

\noindent {\bf Assumptions for Section \ref{section:theo}:} We are considering a sequence of model instances with the following
three assumptions:
\begin{itemize}
\item For each value $\nn =1,2,3,\ldots $, we are considering a model instance with $\nn = \max \{ m,n \}$, so the model parameters $m,n,K,p_{ij},\varrho_{ij}$ (for all relevant $i,j$) are all
functions of $\nn$. (However, in the notation for these parameters, we omit $\nn$ for ease of exposition.)
\item We assume that there exists a real number $q \in (0,\frac{1}{2})$, which is also a function of $\nn$, such that, for all
$\nn =1,2,\ldots$, it holds that  $p_{ij} \in [q,1-q]$  for all $i,j: 1\leq i < j \leq \nn$, and $\varrho_{ij} > (1-2q)^2$ for all  $i,j: 1\leq i < j \leq K$.
 \item
For each $\nn =1,2,3,\ldots $, define $\epsilon := \min_{i,j: 1\leq i < j \leq K} \left ( \varrho_{ij} - (1-2q)^2 \right )$;
of course, $\epsilon$ is a function of $\nn$, and $\epsilon \in (0,1)$. We assume that $\epsilon^2 K \in \omega (\log \nn)$;
specifically, $\lim_{\nn \rightarrow \infty}\frac{\epsilon^2K}{\log \nn}=\infty$.
\end{itemize}

We will say that a property of such random graphs $(G,H)$ holds {\it almost always} to mean that, with probability $1$, the property
holds for all but a finite number of values of $\nn$.

Recall that the matrices $A$ and $B$ are $\pm 1-$adjacency matrices, i.e.~$(1,-1,-1)$-adjacency matrices in the notation of Equation \ref{eqn:adjacency}.
This means that the diagonals of $A$ and $B$ are all $-1$'s. However, as long as the values of the diagonals are all equal to each other, the specific (common) value
does not affect the solution of the subgraph-subgraph matching problem
 \begin{eqnarray}
\max_{X \in \vP_{m,n,K}}\TT AXBX^T \ \
\end{eqnarray}
since it just changes the objective function by a constant. It will be helpful for our
analysis in Section \ref{section:theo} to have the diagonals of $A$ and $B$ be all zeros, and
we adopt this convention. (So $A$ and $B$ are $(1,-1,0)-$adjacency matrices.) The relaxed subgraph-subgraph matching problem
 \begin{eqnarray}
\max_{X \in \vD_{m,n,K}}\TT AXBX^T \ \
\end{eqnarray}
is affected, but the results of this convention do not empirically seem to change the efficacy of the ssSGM algorithm. The main result of
Section \ref{section:theo} is the following. ($I_K$ denotes the $K\times K$ identity matrix, and the zero matrices are appropriately sized.)

\begin{theorem} \label{thm:main}
Under the above Assumptions of Section \ref{section:theo}, it almost always holds that
$$ \arg \max_{X \in \vP_{m,n,K}}\TT AXBX^T = \arg \max_{X \in \vD_{m,n,K}}\TT AXBX^T = \left \{ \II  \right \} . $$
\end{theorem}

In other words, under the above-specified conditions, the subgraph-subgraph matching problem solution almost always locates the correct core and
correctly matches these core vertices. Moreover, even the relaxed problem solution achieves this goal, which further motivates the focus of the ssSGM algorithm on solving the relaxed problem.

The following is the key lemma from which our main result Theorem~\ref{thm:main} quickly follows.
\begin{lemma} \label{lemma:main}
Under the above Assumptions of Section \ref{section:theo}, it almost always holds that
$$ \arg \max_{P,Q \in \vP_{m,n,K}}\TT APBQ^T = \left \{ \II  \right \} . $$
\end{lemma}

Because of its length and complexity, Lemma \ref{lemma:main} is proved in \ref{appendixA}, which the reader is strongly encouraged to read. \\

\noindent {\bf Proof of Theorem \ref{thm:main}:}
By Lemma \ref{lemma:hull}, any $X \in \vD_{m,n,K}$ can be expressed
$X=\sum_i \alpha_i P^{(i)}$ where $\{ \alpha_i \}$ are nonnegative real numbers that sum to $1$ and
$\{ P^{(i)} \} = \vP_{m,n,K}$. Then  $$\TT AXBX^T= \TT A \left ( \sum_i\alpha_iP^{(i)} \right ) B \left ( \sum_j\alpha_jP^{(j)} \right )^T=
\sum_{i,j}\alpha_i \alpha_j \TT A P^{(i)}B P^{(j)T}$$ is a convex combination of expressions of the form
$\TT APBQ^T$ in Lemma \ref{lemma:main}, and Theorem \ref{thm:main} follows immediately from
Lemma \ref{lemma:main}. $\qed$

\section{Discussion, future work \label{sec:discussion}}

The graph matching problem (as well as the many variants and generalizations) is important because of its wide-ranging application as well as its challenging theory and computational complexity.
We provide here an efficient and effective algorithm (ssSGM) to approximately solve the seeded subgraph-subgraph matching problem, and we state and prove an information-theoretic matchability result (Theorem \ref{thm:main})
in this context.

There are a number of promising directions for further work that build on our work here.

The seeds that are incorporated into
ssSGM are vertices in the core together with their a priori known alignments across the pair of graphs; such known alignments are the usual
definition of ``seeds" in graph
matching. However, in our setting
where there are core vertices and non
core vertices, we might know
a priori that certain vertices are
in the core, but not know their alignments.
It would be interesting to
incorporate  such ``pre-seeds" into our
subgraph-subgraph matching problem, and compare their effectiveness to seeds in identifying a latent matching.

Also, Theorem \ref{thm:main} provides
(mild) sufficient conditions for matchability in our setting, but it does not provide necessary conditions.
It would be worth extending our results to obtain necessary and sufficient conditions
(up to order).
In particular, the ``total correlation" (a measurement of correlation between a pair of graphs) introduced in our previous work  \citet{FISHKIND2019295} and \citet{fishkind2021phantom} combines
the edge correlation together with
the Bernoulli coefficients, and
was demonstrated to be an
extremely effective
tool for determining matchability.
The random pair-of-graphs model in those two papers 
was less general than here,
but sets forward a promising tool
to achieve necessary and sufficient
conditions. \\

\noindent
{\bf Acknowledgments:} The authors are grateful to the anonymous referees and editors whose very thoughtful remarks  added much to this work.

\bibliographystyle{asa.bst}
\bibliography{biblio2}

\appendix
\renewcommand{\thesection}{Appendix \Alph{section}}

\section{The Proof of Lemma \ref{lemma:main} \label{appendixA}}

In order to prove Lemma \ref{lemma:main}, which was the fundamental result used just above to prove our main result Theorem~\ref{thm:main}, we provide some preliminaries next, including some observations regarding the distribution of the random graphs and also another key lemma. Then we will prove
Lemma \ref{lemma:main}.\\

\noindent
{\bf Preliminaries for the proof of Lemma \ref{lemma:main}:}\\ 

In this subsection we provide some preliminaries
for the proof of
 Lemma \ref{lemma:main}, continuing with the notation of Section \ref{section:model}.

For any positive integers $i,j$ such that $1\leq i < j \leq K$, it is straightforward to calculate
the joint probabilities of adjacency for $i,j$ in $G$ and $i,j$ in $H$:
\begin{eqnarray}
\begin{array}{c||cc}
    & i \sim_H j   & i \not \sim_H j \\ \hline \hline
  i \sim_G j       & p_{ij}^2+\varrho_{ij}p_{ij}(1-p_{ij}) & (1-\varrho_{ij})p_{ij}(1-p_{ij}) \\
  i \not \sim_G j  & (1-\varrho_{ij})p_{ij}(1-p_{ij})    & (1-p_{ij})^2+\varrho_{ij}p_{ij}(1-p_{ij})
\end{array} \label{eqn:probs}
\end{eqnarray}

Suppose that the random variables $Y \sim \textup{Bernoulli}(p_{ij})$, \ $Y_{Heads} \sim \textup{Bernoulli}\left (\varrho_{ij} \cdot 1 + (1-\varrho_{ij})p_{ij} \right )$, \
$Y_{Tails} \sim \textup{Bernoulli} \left ( \varrho_{ij} \cdot 0 + (1-\varrho_{ij})p_{ij} \right )$ are independent. It will be useful to
note by simple computation that the joint distribution of the pair of random variables $(Y,\ (1-Y)Y_{Tails}+Y \cdot Y_{Heads})$ is that same
as that of $(Y,Y')$ where $Y' \sim \textup{Bernoulli}(p_{ij})$ and the Pearson correlation coefficient for $Y,Y'$ is $\varrho_{ij}$; indeed,
the joint probabilities are given above in (\ref{eqn:probs}). As a consequence, the distribution of the pair of random variables $(A_{ij},B_{ij})$ is a function of three independent Bernoulli random variables. The following result from \citet{kim} will be very useful to us.

\begin{lemma} \label{thm:bound}
  Let $Y_1,Y_2,\ldots,Y_N$ be independent Bernoulli random variables with respective Bernoulli parameters $r_1,r_2,\ldots,r_N$. Let
$g: \{ 0,1\}^N \rightarrow \R$ be a function and $L\in \R_{>0}$ with the property that for any $y,y' \in \{0,1\}^N$ such that $y,y'$ differ in exactly one position (i.e. $\sum_{i=1}^N |y_i-y'_i|=1$) it holds that $|g(y)-g(y')|\leq L$. Define $\sigma:=L \sqrt{\sum_{i=1}^Nr_i(1-r_i)}$ and define
$Y:=g(Y_1,Y_2,\ldots,Y_N)$. Then
$$ \mbox{for any } \tau: 0 \leq \tau < 2\frac{\sigma}{L} \  \mbox{  it holds that  } \ \p \left (  |Y-\e Y | > \tau \sigma   \right )   \leq  2e^{-\frac{\tau^2}{4}}. $$
\end{lemma}
\noindent See \citet{kim} for a proof of Lemma  \ref{thm:bound}.\\

\noindent 
{\bf Proof of Lemma \ref{lemma:main}:}\\  

In this subsection we prove Lemma \ref{lemma:main}, continuing with the notation of Section \ref{section:model}.

For any particular value of $\nn$, consider for now any particular
(fixed) $P,Q \in \vP_{m,n,K}$ such that $P$ and $Q$ are not both $\II$. Define the random variable
\begin{eqnarray}
Y_{PQ}:=  \TT A \II B \II ^T  -  \TT APBQ^T \nonumber
\end{eqnarray}
(which is random by virtue of the $\pm 1$-adjacency matrices  $A$ and $B$); our first task (receiving the bulk of our attention in this proof) is to bound $\p (Y_{PQ} \leq 0)$.
With a strong enough bound, we will then be able to consider all possible such $P,Q$, and all $\nn$, and quickly
prove Lemma \ref{lemma:main}.\\

For any positive integer $a$, recall that $[a]$ denotes the set of integers $\{1,2,3,\ldots,a\}$.  Define
\begin{eqnarray}
\Om_P:= \{ i \in [m] \ : \ \exists j \in [n] \mbox{ such that } P_{ij}=1 \}, \ \  \OM_P:= \{ j \in [n] \ : \ \exists i \in [m] \mbox{ such that } P_{ij}=1 \}, \nonumber \\ \nonumber
\Om_Q:= \{ i \in [m] \ : \ \exists j \in [n] \mbox{ such that } Q_{ij}=1 \},  \ \ \OM_Q:= \{ j \in [n] \ : \ \exists i \in [m] \mbox{ such that } Q_{ij}=1 \} .
\end{eqnarray}
Of course, $|\Om_P|=|\OM_P|=|\Om_Q|=|\OM_Q|=K$, by definition.\\
Define the following bijections:
\begin{eqnarray}
\phi_P:\Om_P \rightarrow \OM_P \mbox{ such that } \phi_P(i)=j \mbox{ precisely when } P_{ij}=1, \nonumber \\
\phi_Q:\Om_Q \rightarrow \OM_Q \mbox{ such that } \phi_Q(i)=j \mbox{ precisely when } Q_{ij}=1, \nonumber
\end{eqnarray}
In what follows, $(i,j)$ will refer to the ordered pair of indices (as opposed to an open interval). Define
\begin{eqnarray}
  S_I &=& \{ (i,j) \in [K]\times [K] \ : \ i \neq j    \}  \nonumber  \\
  S_{PQ} &=& \{  (i,j) \in \Om_P \times \Om_Q \ : \ i \ne j \mbox{ and } \phi_P(i) \ne \phi_Q(j)  \} \nonumber \\
  S &=&  \{ (i,j) \in S_I  \ : \ \phi_P (i) =i \mbox{ and } \phi_Q(j)=j \}  \nonumber \\
  S' &=& \{ (i,j) \in S_I  \ : \ \phi_P (i) =j \mbox{ and } \phi_Q(j)=i \} \nonumber
\end{eqnarray}
Note that $S$ and $S'$ are disjoint, and that $S \cup S' \subseteq S_I$ and $S \cup S' \subseteq S_{PQ}$. Also
note that for all $(i,j) \in S \cup S'$ we have $A_{ij}B_{ij}=A_{ij}B_{\phi_P(i),\phi_Q(j)}$, hence it holds that
\begin{eqnarray}
Y_{PQ} & := & \TT A \II B \II ^T  -  \TT APBQ^T \nonumber \\
   & \equiv & \sum_{(i,j)\in S_I}A_{ij}B_{ij} \ \ \ - \sum_{(i,j) \in S_{PQ}} A_{ij}B_{\phi_P(i),\phi_Q(j)} \nonumber \\
   & = & \sum_{(i,j)\in S_I \backslash (S \cup S') }A_{ij}B_{ij} \ \ \ - \sum_{(i,j) \in S_{PQ}\backslash (S \cup S')  } A_{ij}B_{\phi_P(i),\phi_Q(j)}
   \label{eqn:sum}.
\end{eqnarray}
Next, define
\begin{eqnarray}
  T_{PQ} &=& \{ i \in [K] \ : \ \phi_P(i)=\phi_Q(i)=i \} \nonumber \\
  T_P &=& \{ i \in [K]  \ : \ \phi_P(i)= i \} \ \backslash \ T_{PQ} \nonumber  \\
  T_Q &=& \{ i \in [K] \ : \  \phi_Q(i)= i \} \ \backslash \ T_{PQ} \nonumber
\end{eqnarray}
Note that $T_{PQ} \subseteq \Om_P \cap \Om_Q$, $T_P \subseteq \Om_P$, and $T_Q \subseteq \Om_Q$.
Of course, $T_{PQ}$, $T_P$, and $T_Q$ are mutually disjoint. Denote $t_{PQ}:=|T_{PQ}|$,
$t_P:=|T_P|$, and $t_Q:=|T_Q|$.

The following relationships are easy to demonstrate:
\begin{eqnarray}
  t_{PQ} &<& K \\
  t_{PQ}+t_P+t_Q &\leq & K \label{eqn:four} \\
  |S_I| &=& K^2-K     \label{eqn:nine} \\
  |S_{PQ}| & \leq & K^2-t_{PQ} \label{eqn:one} \\
  |S_{PQ}| & \geq & K^2-2K+t_{PQ} \label{eqn:two} \\
  |S'| & \leq & K-t_{PQ} \\
  |S| &=& t_{PQ}(t_{PQ}-1)+t_{PQ}(t_P+t_Q)+t_Pt_Q \nonumber \\
      & \leq &  t_{PQ}(t_{PQ}-1)+t_{PQ}(K-t_{PQ})+ \frac{(K-t_{PQ})^2}{4} \nonumber \\
      & = & \frac{K^2}{4} +\frac{Kt_{PQ}}{2} + \frac{t_{PQ}^2}{4} -t_{PQ} \label{eqn:three}  
\end{eqnarray}
Inequality (\ref{eqn:one}) holds since $|\Om_P|=|\Om_Q|=K$ and, if $i \in T_{PQ}$, then $(i,i)\in \Om_P \times \Om_Q$ and $(i,i)\not \in S_{PQ}$.
Inequality (\ref{eqn:two}) holds since $\{  (i,j) \in \Om_P \times \Om_Q \ : \ i \ne j  \}$ has at least $K^2-K$ elements, and
$\{  (i,j) \in \Om_P \times \Om_Q \ : \ i \ne j \mbox{ and } \phi_P(i) = \phi_Q(j)  \}$ has at most $K-t_{PQ}$ elements.
The inequality in (\ref{eqn:three}) follows from  (\ref{eqn:four}) and the Arithmetic-Geometric Mean Inequality, which
says $t_Pt_Q \leq \left ( \frac{t_P+t_Q}{2}  \right )^2$.

Now, define $\beta$ and $\gamma$ as follows:
\begin{eqnarray}
  \beta &:=& K^2-3K+2t_{PQ}- (  t_{PQ}(t_{PQ}-1)+t_{PQ}(t_P+t_Q)+t_Pt_Q   )  \nonumber \\
  \gamma &:=&  K^2-t_{PQ} -  (  t_{PQ}(t_{PQ}-1)+t_{PQ}(t_P+t_Q)+t_Pt_Q   ) \label{eqn:beta}
\end{eqnarray}
The point of defining $\beta$ and $\gamma$ as such is that, by easy application of Equations (\ref{eqn:nine}) through (\ref{eqn:three}),
the cardinalities of the sets indexing the sums in Equation \ref{eqn:sum} can be bounded as
\begin{eqnarray}
  \beta & \leq & |S_I \backslash (S \cup S')| = |S_I|-|S|-|S'| \leq   \gamma  \label{eqn:seven} \\
  \beta & \leq & |S_{PQ} \backslash (S \cup S')| = |S_{PQ}|-|S|-|S'| \leq   \gamma  \label{eqn:eight} .
\end{eqnarray}

We next show that, for $\nn$ large enough, it holds that
\begin{eqnarray}
\beta > 0 \ \ \ \mbox{ and } \ \ \ \gamma \leq 2 \beta .  \label{eqn:gamma}
\end{eqnarray}
Indeed, using the definition of $\beta$ in (\ref{eqn:beta}), applying the inequality (\ref{eqn:three}), gives us
\begin{eqnarray}
  \beta
   & \geq &  K^2-3K+2t_{PQ}- \left (  \frac{K^2}{4} +\frac{Kt_{PQ}}{2} + \frac{t_{PQ}^2}{4} -t_{PQ} \right  ) \nonumber  \\
   & = & \frac{3}{4}K^2-3K+3t_{PQ}-\frac{Kt_{PQ}}{2}-\frac{t_{PQ}^2}{4}   \nonumber  \\
   & \geq & \frac{3}{4}K^2-3K+3t_{PQ}-\frac{Kt_{PQ}}{2}-\frac{K^2}{4} \nonumber \\
   & = & \frac{1}{2} (K-t_{PQ})(K-6) \  > 0 \ \ \ \ \ \ \ \ \mbox{ when } K > 6.   \label{eqn:ten}
\end{eqnarray}
From the Assumptions of Section \ref{section:theo} listed in Section \ref{section:model}, recall that
$\epsilon^2 K \in \omega (\log \nn)$. Since $\epsilon \in (0,1)$, we have that $K$ goes to infinity as $\nn$ goes to infinity.
Thus $\beta > 0$ when $\nn$ is large enough. Next, using the definitions of $\beta$ and $\gamma$ in (\ref{eqn:beta}), applying the inequality (\ref{eqn:three}), gives us
\begin{eqnarray}
  2 \beta - \gamma &=& K^2-6K+5t_{PQ}- (t_{PQ}(t_{PQ}-1)+t_{PQ}(t_P+t_Q)+t_Pt_Q)  \nonumber \\
    &\geq & K^2-6K+5t_{PQ}- \left ( \frac{K^2}{4} +\frac{Kt_{PQ}}{2} + \frac{t_{PQ}^2}{4} -t_{PQ} \right ) \nonumber  \\
    &=& \frac{3}{4}K^2 -6K+6t_{PQ}-\frac{Kt_{PQ}}{2}-\frac{t_{PQ}^2}{4}. \label{eqn:bottom}
\end{eqnarray}
The bottom line of (\ref{eqn:bottom}) is a concave quadratic in the variable $t_{PQ}$ hence, when minimizing it over $t_{PQ}$ in the range of
$t_{PQ}\in \{ 0,1,2,\ldots,K-1\}$, the minimum is achieved when $t_{PQ}$ is $0$ or $K-1$. Plugging in $0$ for $t_{PQ}$ in the
bottom line of (\ref{eqn:bottom}) yields $\frac{3}{4}K^2-6K$, and plugging in $K-1$ for $t_{PQ}$ in the
bottom line of (\ref{eqn:bottom}) yields $K-6-\frac{1}{4}$. Hence we
have that $2 \beta -\gamma \geq 0$ for $\nn$ big enough (hence $K$ is big enough), completing our demonstration of (\ref{eqn:gamma}).

We are now almost prepared to bound the expectation of $Y_{PQ}$.
From the joint probability formula in Equation(\ref{eqn:probs}), we calculate, for
each $(i,j) \in S_I$,
\begin{eqnarray}
\e (A_{ij}B_{ij}) &=&  p^2_{ij}+\varrho_{ij}p_{ij}(1-p_{ij}) + (1-p_{ij})^2+\varrho_{ij}p_{ij}(1-p_{ij})
- 2   (1-\varrho_{ij})p_{ij}(1-p_{ij}) \nonumber  \\
 &=& (1-2p_{ij})^2+4 p_{ij}(1-p_{ij})\varrho_{ij} \nonumber \\
  &\geq & (1-2p_{ij})^2\varrho_{ij}+4 p_{ij}(1-p_{ij})\varrho_{ij} \nonumber \\
 &=& \varrho_{ij} \ \ \geq (1-2q)^2 + \epsilon , \label{eqn:six}
\end{eqnarray}
the last inequality following from the Assumptions for Section \ref{section:theo} listed in Section \ref{section:model}.
However, for any positive integers $i,j \in [m]$ such that $i \ne j$ and any
$i',j'\in [n]$ such that $i' \ne j'$, and also such that $i$ being adjacent to $j$ in $G$
is independent of $i'$ being adjacent to $j'$ in $H$, we have
\begin{eqnarray} \label{eqn:five}
  \e (A_{ij}B_{i'j'}) \ = \ \e A_{ij} \ \e B_{i'j'} \ = \ (1-2p_{ij})(1-2p_{i'j'}) \ \leq (1-2q)^2 ,
\end{eqnarray}
the last inequality following from the Assumptions for Section \ref{section:theo} listed  in Section \ref{section:model}.
Using Equations  (\ref{eqn:six}), (\ref{eqn:five}), (\ref{eqn:seven}), and (\ref{eqn:eight}), we take the expectation of $Y_{PQ}$ in Equation (\ref{eqn:sum}) and bound it
as follows.
\begin{eqnarray}
  \e Y_{PQ} &=&   \sum_{(i,j)\in S_I \backslash (S \cup S') }\e A_{ij}B_{ij}\ \ \ - \sum_{(i,j) \in S_{PQ}\backslash (S \cup S')  } \e A_{ij}B_{\phi_P(i),\phi_Q(j)} \nonumber \\
   & \geq & \sum_{(i,j)\in S_I \backslash (S \cup S') } \left ( (1-2q)^2+\epsilon \right ) \ \ \   - \sum_{(i,j) \in S_{PQ}\backslash (S \cup S')  } (1-2q)^2 \nonumber \\
   &\geq & \left ( (1-2q)^2+ \epsilon \right ) \beta - (1-2q)^2 \gamma \nonumber \\
   &=& (1-2q)^2(\beta-\gamma)+\epsilon \beta \nonumber \\
   &=& \epsilon \beta -3(1-2q)^2(K-t_{PQ}). \label{eqn:expect}
\end{eqnarray}
From (\ref{eqn:expect}), and the inequality (\ref{eqn:ten}) for $\beta$, we obtain that
\begin{eqnarray}
  \e Y_{PQ} - \frac{\epsilon \beta}{2} & \geq & \frac{\epsilon \beta}{2} - 3(K-t_{PQ})  \nonumber \\
   & \geq & \frac{\epsilon}{2} \left ( \frac{3}{4}K^2-3K+3t_{PQ}-\frac{Kt_{PQ}}{2}-\frac{t_{PQ}^2}{4} \right )    -3(K-t_{PQ}) \nonumber \\
   & \geq & \frac{3}{8}\epsilon K^2 -6(K-t_{PQ})-\frac{1}{4}\epsilon K t_{PQ} - \frac{1}{8}\epsilon t_{PQ}^2 \label{eqn:eleven}
\end{eqnarray}
The bottom line of (\ref{eqn:eleven}) is a concave quadratic in the variable $t_{PQ}$ hence, when minimizing it over $t_{PQ}$ in the range of
$t_{PQ}\in \{ 0,1,2,\ldots,K-1\}$, the minimum is achieved when $t_{PQ}$ is $0$ or $K-1$. Plugging in $0$ for $t_{PQ}$ in the
bottom line of (\ref{eqn:eleven}) yields $\frac{3}{8}\epsilon K^2-6K=K(\frac{3}{8}\epsilon K-6) $, and plugging in $K-1$ for $t_{PQ}$ in the
bottom line of (\ref{eqn:eleven}) yields an expression bounded below by $\frac{1}{2}\epsilon K -7$.
However,  $\epsilon^2 K \in \omega (\log \nn)$ and  $\epsilon \in (0,1)$ imply that $\epsilon K$ is greater than any particular number when $\nn$ is large  enough, hence for large enough $\nn$ we have that the expression in (\ref{eqn:eleven}) is positive, i.e.
\begin{eqnarray}
  \e Y_{PQ}  > \frac{\epsilon \beta}{2}. \label{eqn:fifteen}
\end{eqnarray}

We are now in a position to bound $\p (Y_{PQ} \leq 0)$, which will be accomplished using Lemma \ref{thm:bound}; there are a number of
conditions and requirements for Lemma \ref{thm:bound}, which we next verify as holding.

 For each
$i,j: 1 \leq i <j \leq K$, recall from the paragraph immediately preceding Lemma \ref{thm:bound} that the distribution of
the pair of random variables $(A_{ij},B_{ij})$ may be thought of as
 a function of three independent Bernoulli random variables, one of them having
Bernoulli parameter in the interval $[q,1-q]$. Indeed, the joint distribution of all $A_{ij}, B_{i'j'}$, collectively for all
$i,j: 1 \leq i < j \leq m$ and $i',j':1 \leq i'<j'\leq n$ is thus a function of $3\binom{K}{2}+ \left [  \binom{m}{2} - \binom{K}{2}  \right ]
+\left [   \binom{n}{2} - \binom{K}{2}   \right ]$ collectively independent Bernoulli random variables.

Specifically, by (\ref{eqn:sum}), we have that $Y_{PQ}$ is a function of such independent Bernoulli random variables which underlie $A_{ij}B_{ij}$
for all  $(i,j)\in S_I \backslash (S \cup S')$ and which underlie  $A_{ij}B_{\phi_P(i),\phi_Q(j)}$ for all $(i,j) \in S_{PQ}\backslash (S \cup S')$.  Let us denote by $N$ the number of these Bernoulli random variables; from (\ref{eqn:seven}), (\ref{eqn:eight}), and (\ref{eqn:gamma}) we have that
\begin{eqnarray}
 \frac{3}{2} \beta \leq N \leq 6 \gamma \leq 12 \beta  \label{eqn:thirteen}
\end{eqnarray}
(since we are allowing for three such Bernoullis for each such $(i,j)$, and there is a $\frac{1}{2}$ in the left-hand side because the $(i,j)$ are {\it ordered} pairs).
Denote by $r_1,r_2,\ldots,r_N$ the associated  Bernoulli parameters; at least $\frac{\beta}{2} $ of these are in the interval $[q,1-q]$. Thus
\begin{eqnarray}
 q(1-q) \frac{\beta}{2} \leq \sum_{i=1}^Nr_i(1-r_i) \leq \sum_{i=1}^N \frac{1}{2}(1-\frac{1}{2})=\frac{N}{4} . \label{eqn:fourteen}
\end{eqnarray}
It is easy to see that changing the realized value of any of these $N$ Bernoulli random variables (from zero to one or vice versa) can change the realized value of $ Y_{PQ}$ by at most $L:=8$; define $\sigma:=L \sqrt{\sum_{i=1}^Nr_i(1-r_i)}$, hence $\sigma^2=64\sum_{i=1}^Nr_i(1-r_i)$.
 Thus, by (\ref{eqn:thirteen}) and (\ref{eqn:fourteen}),\  $\sigma^2$ satisfies
\begin{eqnarray}
  32 q(1-q) \beta \ \leq \  \sigma^2 \ \leq \ 192 \beta .  \label{eqn:twelve}
\end{eqnarray}
Finally, since by definition of $\epsilon$ we have $0< \epsilon \leq 1-(1-2q)^2$, it follows from (\ref{eqn:twelve}) and (\ref{eqn:gamma}) that
\begin{eqnarray}
  \frac{1}{2}\sigma^2 -\epsilon \beta & \geq & 16 q(1-q)\beta - \left [ 1-(1-2q)^2 \right ] \beta \nonumber \\
        & = & 12 q(1-q)\beta >0, \nonumber
\end{eqnarray}
from which it immediately follows that (since $L=8$)
\begin{eqnarray}
  \mbox{if we define } \tau:=\frac{\epsilon \beta}{2 \sigma} \ \mbox{ then } \ \tau < 2\frac{\sigma}{L}.
\end{eqnarray}
We now apply Lemma \ref{thm:bound} with $\tau:=\frac{\epsilon \beta}{2 \sigma}$ and $\sigma$ as above, since the lemma's conditions have just been demonstrated to apply; there exists a positive constant $C_1$ (not dependent on $\nn$) such that
\begin{eqnarray}
  \p (Y_{PQ} \leq 0) & \leq & \p \left ( |Y_{PQ}- \e Y_{PQ}| > \frac{\epsilon \beta}{2} \right )  \ \ \ \  \ \  \ \ \mbox{by Equation (\ref{eqn:fifteen})}   \nonumber \\
  &\leq&  2e^{-\frac{\epsilon^2 \beta^2}{16 \sigma^2}} \ \ \ \ \ \ \  \ \ \ \ \ \ \ \ \ \ \ \ \ \ \ \ \ \ \ \ \ \ \ \ \ \mbox{by Lemma \ref{thm:bound}} \nonumber \\
  &\leq&  2e^{-2C_1\epsilon^2 \beta} \ \ \ \ \ \ \  \ \ \ \ \ \ \ \ \ \ \ \ \ \ \ \ \ \ \ \ \ \ \mbox{by Equation (\ref{eqn:twelve}}) \nonumber \\
   &\leq&  2e^{-C_1\epsilon^2 (K-t_{PQ})(K-6)}  \ \ \ \ \ \ \ \ \ \ \ \ \ \  \ \ \mbox{by Equation (\ref{eqn:ten}}) \label{eqn:eighteen}
\end{eqnarray}
when $K>6$.

Until this point, $P$ and $Q$ have been fixed. Keep $\nn$ fixed. For each $t=0,1,2,\ldots,K-1$, define  $\EE_{t,\nn}$ to be
the event that there exists an ordered pair $(P,Q)$ in the set $\Phi^{(t)}$ such that $Y_{PQ}\leq 0$, where
\begin{eqnarray}
\Phi^{(t)} :=  \{ \  (P',Q') \in  \vP_{m,n,K} \times \vP_{m,n,K} \ : \ |\{ i \in [K]  : \phi_{P'}(i)=\phi_{Q'}(i)=i    \}|=t  \ \}
\nonumber
\end{eqnarray}
For each $t=0,1,2,\ldots,K-1$, we have that
\begin{eqnarray}
| \Phi^{(t)}| &\leq & \binom{K}{t} \left [  \binom{m}{K-t} \binom{n}{K-t} (K-t)!  \right ] ^2 \nonumber  \\
& = & \frac{K!}{t!(K-t)!)} \left [   \frac{m!}{(K-t)!(m-(K-t))!} \ \  \frac{n!}{(K-t)!(n-(K-t))!} \ (K-t)!    \right ] ^2 \nonumber \\
& \leq & K^{K-t} \left [  m^{K-t} n^{K-t}  \right ]^2 \nonumber \\
& \leq &  \nn^{5(K-t)} .  \label{eqn:seventeen}
\end{eqnarray}

Applying (\ref{eqn:eighteen}) and (\ref{eqn:seventeen}), there exists a positive constant $C_2$ not dependent on $\nn$ such that when $\nn$ is big enough
it holds that
\begin{eqnarray}
  \p [ \EE_{t,\nn}] & \leq & 2  \ \nn^{5(K-t)} \   e^{-C_1\epsilon^2 (K-t)(K-6)} \ \ = \ \  2e^{5(K-t)\log \nn -C_1\epsilon^2 (K-t)(K-6)} \nonumber \\
  & = & 2 e^{(K-t) \left [  5 \log \nn - C_1 \epsilon^2 (K-6)  \right ]} \nonumber \\
  & \leq & 2 e^{(K-t)(-C_2 \epsilon^2 K)}   \ \ \ \ \ \ \ \ \ \ \ \ \ \ \ \ \ \ \ \ \ \ \ \ \ \ \ \mbox{since } \epsilon^2K \in \omega (\log \nn)  \nonumber \\
  & \leq & 2 e^{(-C_2 \epsilon^2 K)}  \ \ \ \ \ \ \ \ \ \ \ \ \ \ \ \ \ \ \ \ \ \ \ \ \ \ \ \ \ \ \ \ \ \ \mbox{since } K-t \geq 1 \nonumber \\
  & \leq & \frac{1}{\nn^3}   \ \ \ \ \ \ \ \ \ \ \ \ \ \ \ \ \ \ \ \ \ \ \ \ \ \ \ \ \ \ \  \ \ \ \ \ \ \ \ \ \ \ \ \  \mbox{since } \epsilon^2K \in \omega (\log \nn).   \label{eqn:nineteen}
\end{eqnarray}
Hence, we have from (\ref{eqn:nineteen}) that
\begin{eqnarray}
  \p \left [  \bigcup_{t=0}^{K-1}\EE_{t,\nn} \right ] \  \leq K \ \frac{1}{\nn^3} \ \leq \ \frac{1}{\nn^2}; \nonumber
\end{eqnarray}
thus, since $\sum_{\nn=1}^\infty \frac{1}{\nn^2} < \infty $,  Lemma \ref{lemma:main} follows from applying the Borel-Cantelli Lemma to the set of events $\{  \bigcup_{t=0}^{K-1}\EE_{t,\nn}  \}_{\nn=1}^{\infty}$. $\qed$

\section{Proof of  Lemma \ref{lemma:assign}  }
\label{app:pfs}

\noindent {\bf Proof of Lemma \ref{lemma:assign}:}
Suppose $c,d,e$ are any positive integers such that $e\leq \min \{c,d\}$, and consider any particular $M \in \R^{c \times d}$.
Let $\lb$ be any real number less than all of the entries of $M$, and let $\ub$ be any real number greater than
all of the entries of $M$. Define the matrix $\mathcal{M} \in \R^{(c+d-e)\times (c+d-e)}$ such that
$\mathcal{M}= \bigl [\begin{smallmatrix} M & \ub E_{c \times (c-e)} \\ \ub E_{(d-e)\times d} & \ \ \lb E_{(d-e)\times (c-e)} \end{smallmatrix} \bigr ]$, where $E$ is a matrix of all $1$s, sized per subscript.\\
 Consider the linear assignment problem $\max_{\XX \in \vP_{c+d-e,c+d-e,c+d-e}}\TT \mathcal{M}^T \XX$, and also any solution $\mathcal{X} \in \vP_{c+d-e,c+d-e,c+d-e}$, say it is partitioned
as $\XX= \bigl [\begin{smallmatrix} \XX^{(11)} & \XX^{(12)} \\ \XX^{(21)} & \XX^{(22)} \end{smallmatrix} \bigr ]$,
where $\XX^{(11)}, \XX^{(12)}, \XX^{(21)}, \XX^{(22)}$ are respectively of sizes
$c \times d $, and $c \times (c-e)$, and $(d-e)\times d$, and $(d-e)\times (c-e)$. If $\ub$ is great enough,
then it is clear that $\XX^{(12)}$ must have a $1$ in {\bf each} of its $c-e$ columns, and in $c-e$ of its $c$ rows, with the same constant
contribution of $\ub (c-e)$ to the objective function value without regard for which $e$ of its rows do not have $1$s.
Similarly, $\XX^{(21)}$ must have a $1$ in {\bf each} of its $d-e$ rows, and in $d-e$ of its $d$ columns, with the same constant
contribution of $\ub (d-e)$ to the objective function value without regard for which $e$ of its columns do not have $1$s. Thus
$\XX^{(22)}$ is all zeros, and $\XX^{(11)}$ is in $\vP_{c,d,e}$. Ignoring a constant translation of the objective function value
by $\ub(c-e+d-e)$, we thus have that the optimization problem is equivalent to the generalized linear assignment problem
$\max_{X \in \vP_{c,d,e}} \TT M^TX$, whose solution is $\XX^{(11)}$ in the solution $\XX$ of the former problem.\\
Now, simple algebra demonstrates that the above analysis applies equally well as long as $\ub$ is greater than all of the
values of $M$, and the result is shown. $\qed$

\section{More experiments in random graph models}  \label{appendixC}

In this appendix we perform more experiments to showcase the 
effectiveness of ssSGM, in the 
setting of stochastic block 
models and in the setting 
of correlated Bernoulli random graph models.

We begin with 
a stochastic block model.
Say graphs $G$ and $H$
have respective vertex sets $V$ and $W$ 
where $V$ is partitioned into three blocks
$V^{(1)}=\{v_1,v_2,\ldots,v_{100} \}$, 
$V^{(2)}=\{v_{101}, v_{102},\ldots,v_{300} \}$, $V^{(3)}=\{v_{301},
v_{302},\ldots,v_{500} \}$, 
and $W$ is partitioned into three 
blocks $W^{(1)}=\{w_1 ,w_2,\ldots,w_{100}\}$, 
$W^{(2)} 
= \{w_{101}, w_{102},\ldots,w_{300} \}$,
$W^{(3)}=\{w_{301}, w_{302},\ldots,w_{500} \}$.
For any distinct vertices $v,v'\in V$, say $v \in V^{(i)}$ and 
$v' \in V^{(i')}$, the probability that $v$ is adjacent to $v'$ in $G$ is $P_{i,i'}$ for 
\begin{align*}
    P = \begin{bmatrix}
0.5  & 0.3 & 0.3\\
0.3 & 0.45  & 0.3\\
0.3 & 0.3 & 0.4
\end{bmatrix}
\end{align*}
and, for any distinct vertices $w,w'\in W$, say $w \in W^{(i)}$ and 
$w' \in W^{(i')}$, the probability that $w$ is adjacent to $w'$ in $H$ is $P_{i,i'}$. Let 
$\mathcal{I}:=[100]$; 
for any pair of distinct 
$j,j' \in \mathcal{I}$, 
the Pearson correlation coefficient for the Bernoulli random variable 
indicators ${\bf 1}_{v_j \sim_G v_{j'}}$ and 
${\bf 1}_{w_j \sim_H w_{j'}}$ is  $\varrho = .75$, and all other adjacencies are independent. Thus, the core
is $\{ v_j: j \in \mathcal{I} \}$ and $\{ w_j: j \in \mathcal{I} \}$ with the association $v_j \leftrightarrow w_j$.
For each $s=5,10,15,\ldots,50$,
we did 200 experiments of
instantiating random such ($G$,$H$) randomly, selecting $s$ seeds from the core, and running the ssSGM Algorithm and the SGM Algorithm; we plot the
average match ratio vs $s$ in the left panel of Figure 
\ref{fig:sbm}.
\begin{figure}[h!]
    \centering
    \includegraphics[width=0.45\linewidth]{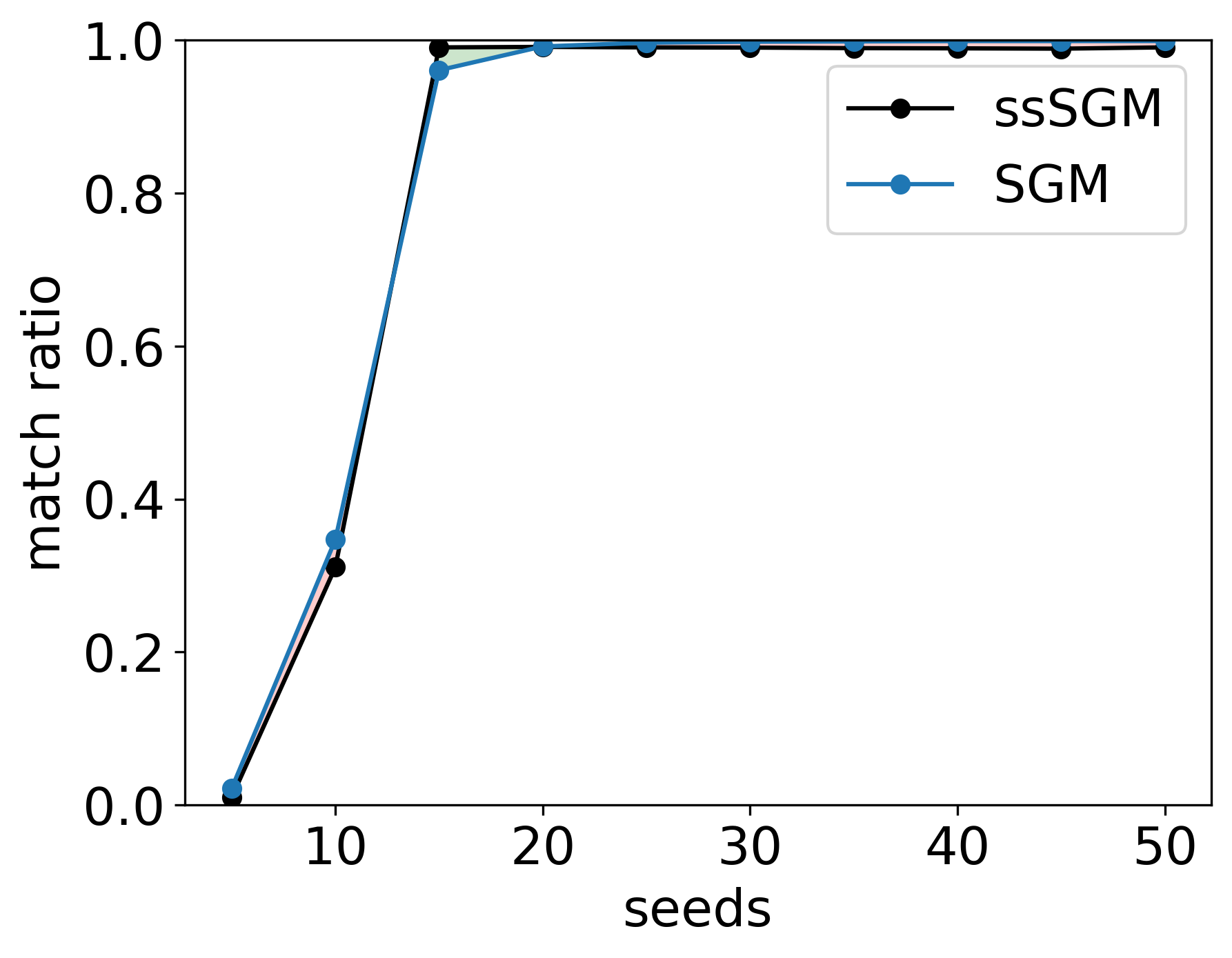}
    \includegraphics[width=0.45\linewidth]{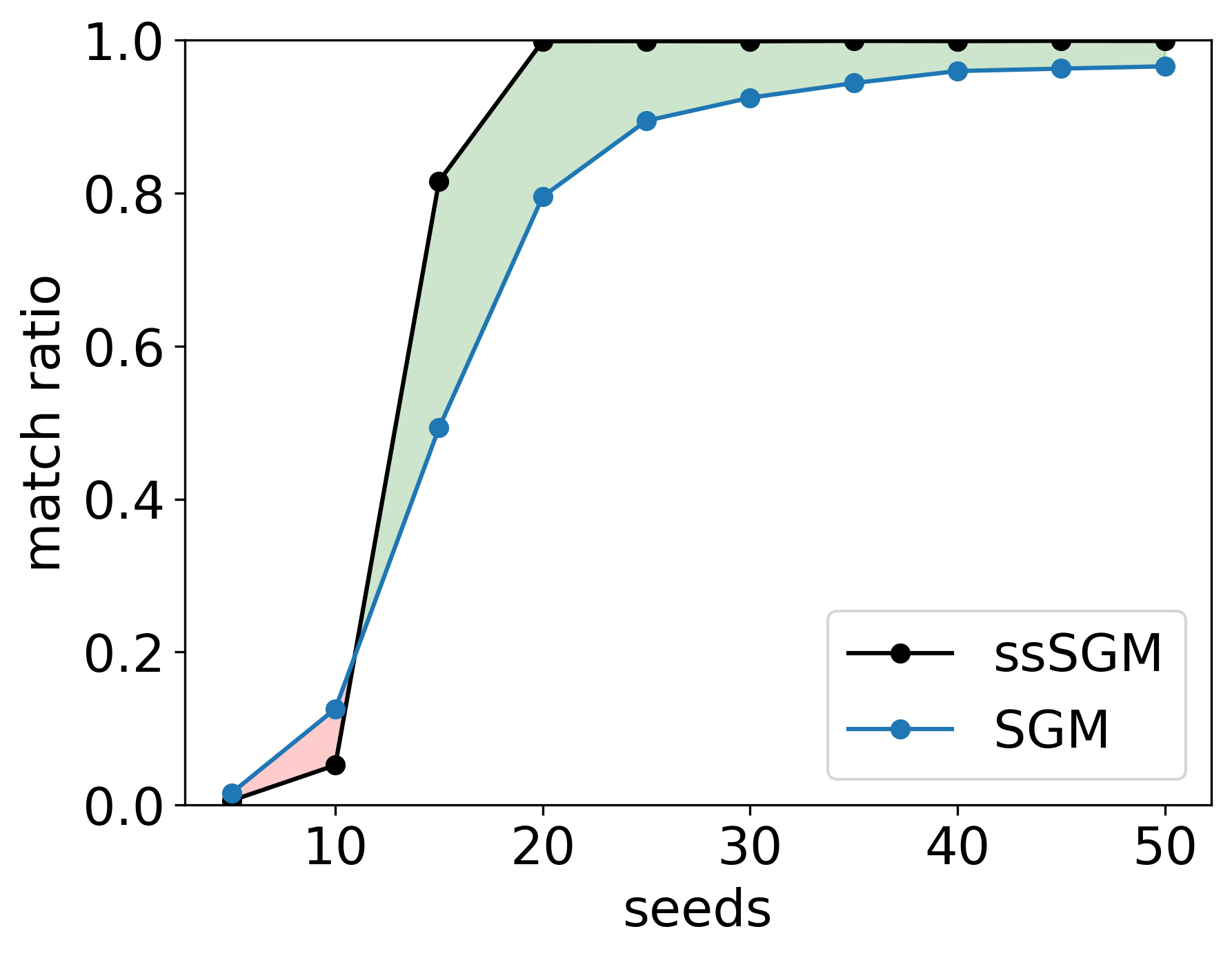}
    \caption{Stochastic block 
    model experiments with core being one of the blocks (left panel) and from the 
    different blocks (right panel).}
    \label{fig:sbm}
\end{figure}

We then repeated the above experiments, the only difference now being that we let
$\mathcal{I}= \{1,2,\ldots,50,\ \ 101,102,\ldots,125, \ \ 301,302,\ldots,325 \}$; 
we plot the
average match ratio vs $s$ in the right panel of Figure 
\ref{fig:sbm}. In both settings 
the core is $100$ vertices, but in the former the core is one block and in the latter the core is from the different blocks.\\

Next we conducted experiments in  correlated Bernoulli random graph models.
Graphs $G$ and $H$ have respective vertex sets  
$V=\{v_1,v_2,\ldots,v_{500}\}$ and
$W=\{w_1,w_2,\ldots,w_{500}\}$. For all integers $i,j$ such that $1 \leq i <j \leq 500$,
we sample $p_{ij}$ from a Uniform$(0.1,0.3)$ distribution, 
and for all  integers $i,j$ such that $1 \leq i <j \leq 100$,
we sample $\varrho_{ij}$ from a Uniform$(0.7,0.8)$ distribution.
Then, for all integers $i,j$ such that $1 \leq i <j \leq 500$, 
the probability that $v_i \sim v_j$ in $G$ is $p_{ij}$, the probability 
that $w_i \sim w_j$ in $H$ is $p_{ij}$ and, for all integers 
$i,j$ such that $1 \leq i <j \leq 100$, 
the Pearson correlation coefficient for the Bernoulli random variable 
indicators ${\bf 1}_{v_j \sim_G v_{j'}}$ and 
${\bf 1}_{w_j \sim_H w_{j'}}$ is  $\varrho_{ij}$, all other adjacencies are 
independent. Thus, the core is $\{v_1,v_2,\ldots,v_{100}\}$ and 
$\{w_1,w_2,\ldots,w_{100}\}$.

For each $s=5,10,15, \ldots, 50$, we did $200$ instantiations of $G$ and $H$, randomly choosing $s$ seeds from the core, and ran ssSGM and SGM with $G$ and $H$; in the left 
panel of Figure \ref{fig:cor-ber} we plotted the average match ratio for ssSGM and SGM vs
the number of seeds $s$.

We repeated the above, the only difference being that $p_{ij}$ were distributed 
Uniform$(0.3,0.5)$, and in the middle panel of Figure \ref{fig:cor-ber} we plotted the average match ratio for ssSGM and SGM vs the number of seeds $s$. 
We repeated the above, the only difference being that $p_{ij}$ were distributed 
Uniform$(0.1,0.5)$, and in the right panel of Figure \ref{fig:cor-ber} we plotted the average match ratio for ssSGM and SGM vs the number of seeds $s$.
\begin{figure}[ht]
    \centering
    \begin{subfigure}[b]{0.32\textwidth}
    \centering
    \includegraphics[width=\textwidth]{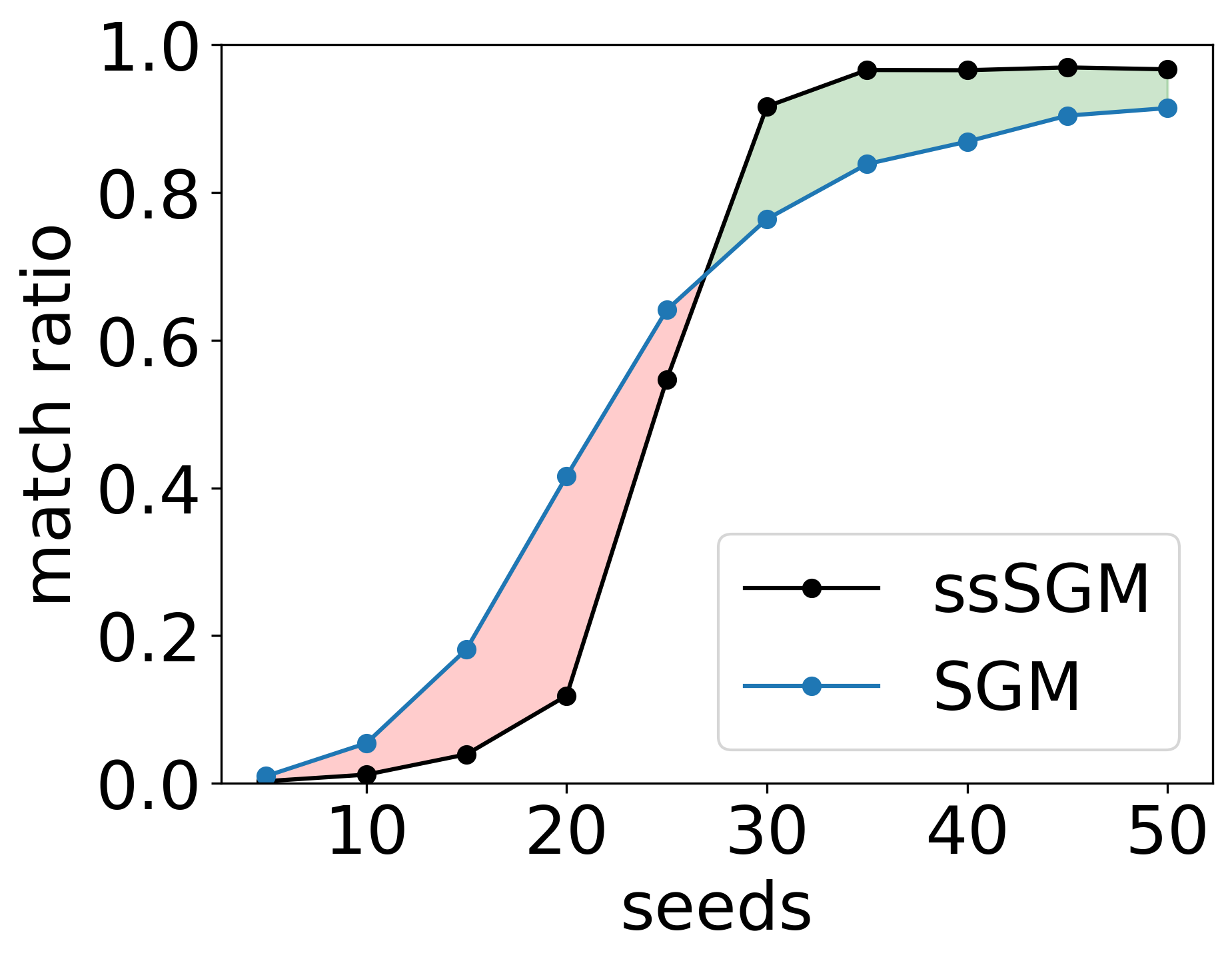}
    \caption*{$p_{ij}\sim \text{Uniform(0.1, 0.3)}$}
    \end{subfigure}
    \hfill
    \begin{subfigure}[b]{0.32\textwidth}
    \centering
    \includegraphics[width=\textwidth]{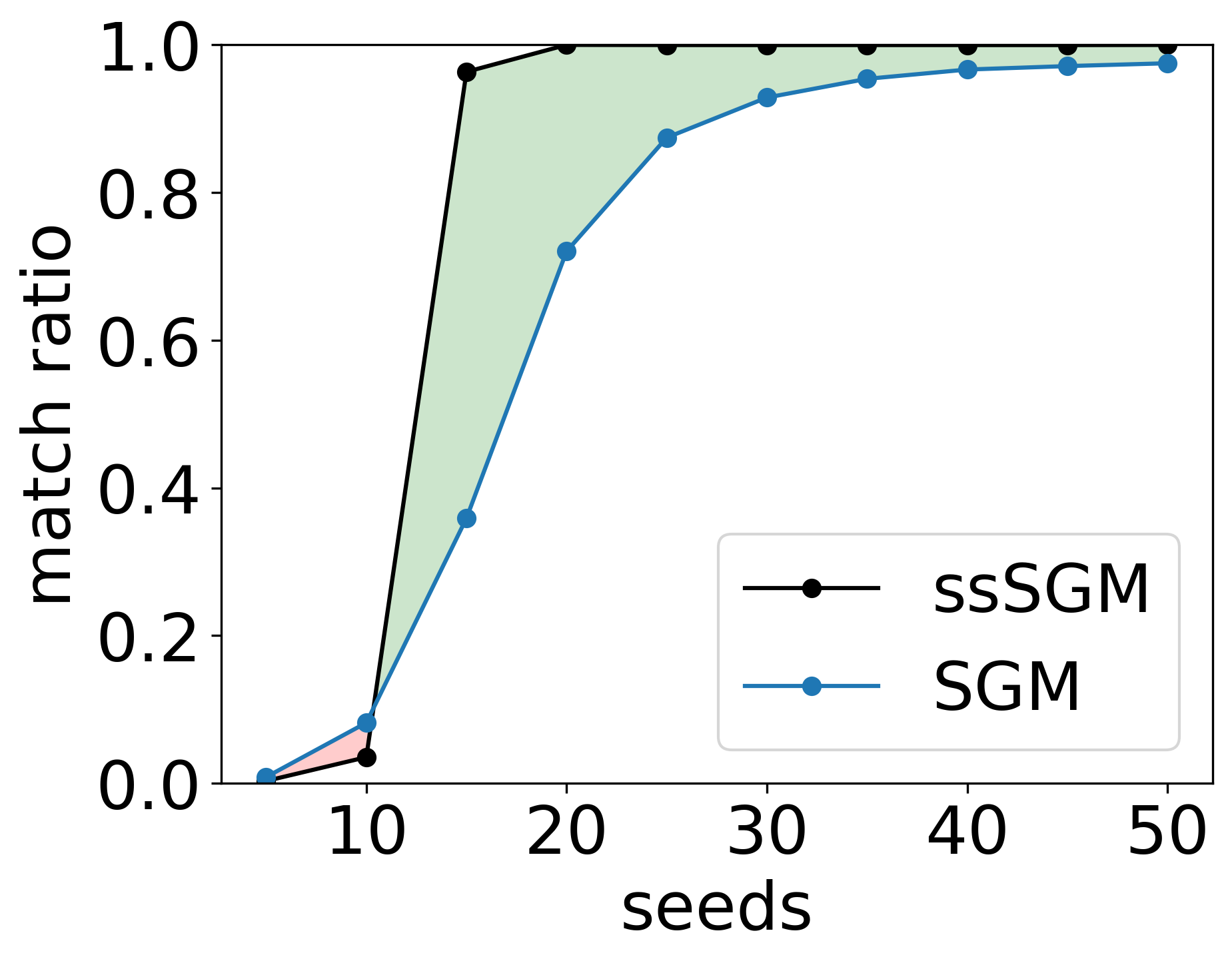}
    \caption*{$p_{ij}\sim \text{Uniform(0.3, 0.5)}$}
    \end{subfigure}
    \hfill
    \begin{subfigure}[b]{0.32\textwidth}
    \centering
    \includegraphics[width=\textwidth]{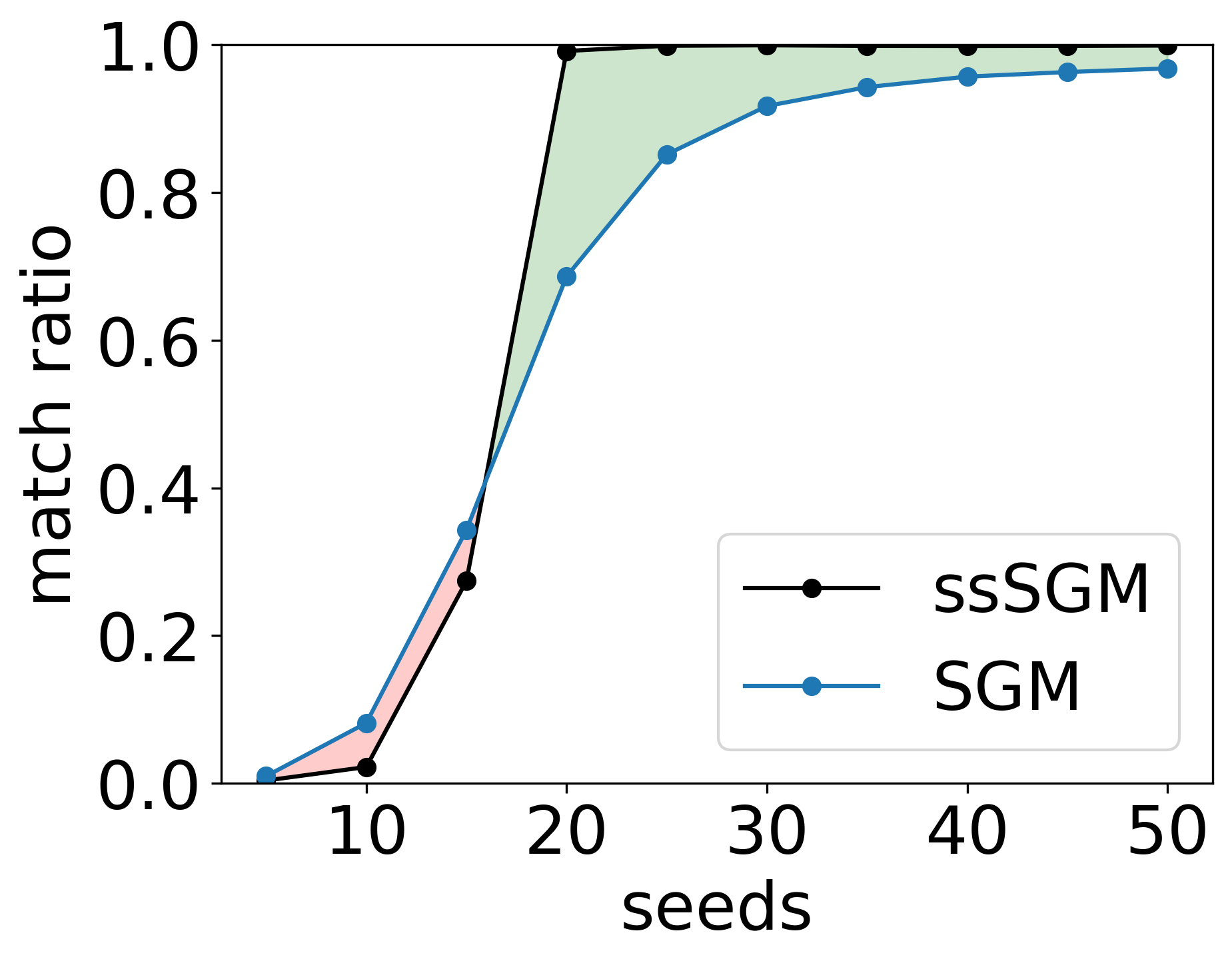}
    \caption*{$p_{ij}\sim \text{Uniform(0.1, 0.5)}$}
    \end{subfigure}
    \caption{Experiments comparing ssSGM and SGM on correlated Bernoulli random graphs}
    \label{fig:cor-ber}
\end{figure}

\section{Comparing use of
$\pm 1$-adjacency matrices vs $\{0,1\}$-adjacency when performing SGM}  \label{appendixD}
In this appendix we provide a brief empirical demonstration to support our claim that the effectiveness of SGM 
is not altered if 
$\pm 1$-adjacency matrices represent the two graphs in the SGM Algorithm vs if $\{ 0,1\}$-adjacency matrices represent the two graphs.

We repeat the experiments from the beginning of Section \ref{sec:exp}, which were illustrated in Figure \ref{fig:synth}, the only difference is that we ran SGM
using $\pm 1$-adjacency matrices to represent the graphs and again with $\{ 0,1\}$-adjacency matrices to represent the graphs, and 
plotted in Figure \ref{fig:SGM_adj} the average 
match ratio as a function of the number of seeds for both 
types of adjacency matrices in SGM;
the left panel comes from using
$p = 0.1, \varrho = 0.5$,
the center panel comes from 
using $p = 0.3, \varrho = 0.6$, and the right panel comes from using $p = 0.5, \varrho = 0.7$. Indeed, there is hardly a difference between the performance of SGM from one type of adjacency matrices to the other.
\begin{figure}[ht]
    \centering
    \begin{subfigure}[b]{0.32\textwidth}
    \centering
    \includegraphics[width=\textwidth]{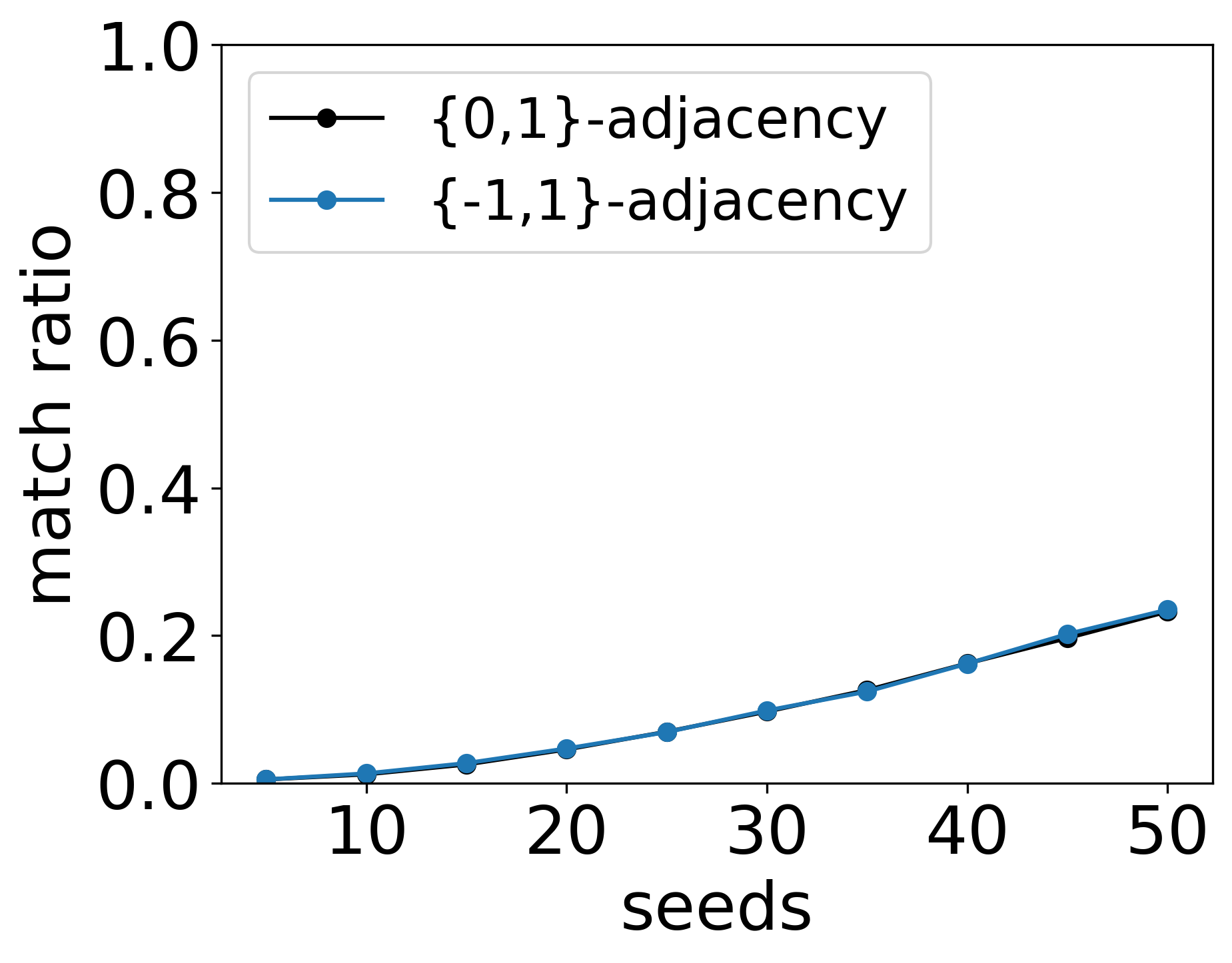}
    \caption*{$p = 0.1, \varrho = 0.5$}
    \end{subfigure}
    \hfill
    \begin{subfigure}[b]{0.32\textwidth}
    \centering
    \includegraphics[width=\textwidth]{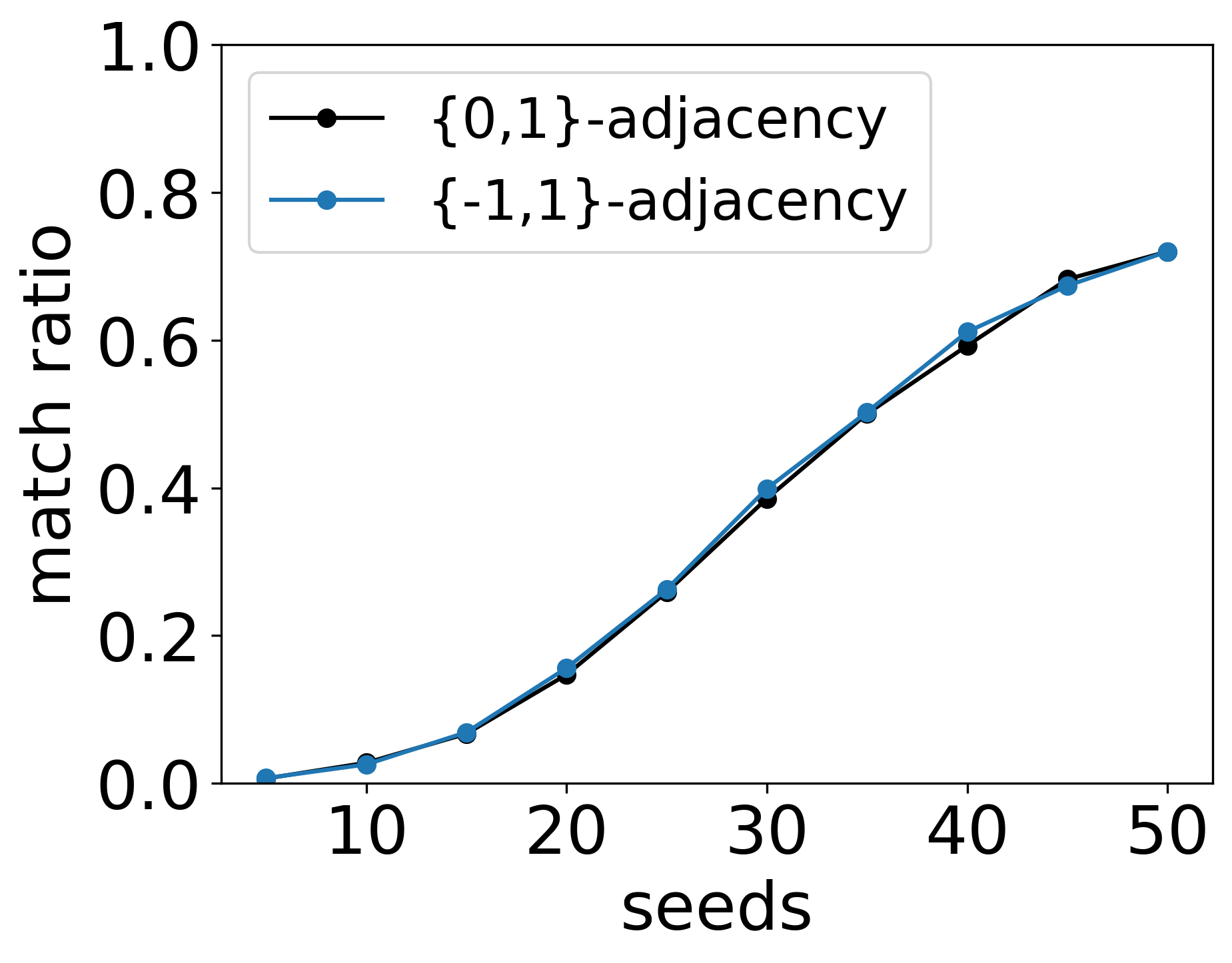}
    \caption*{$p = 0.3, \varrho = 0.6$}
    \end{subfigure}
    \hfill
    \begin{subfigure}[b]{0.32\textwidth}
    \centering
    \includegraphics[width=\textwidth]{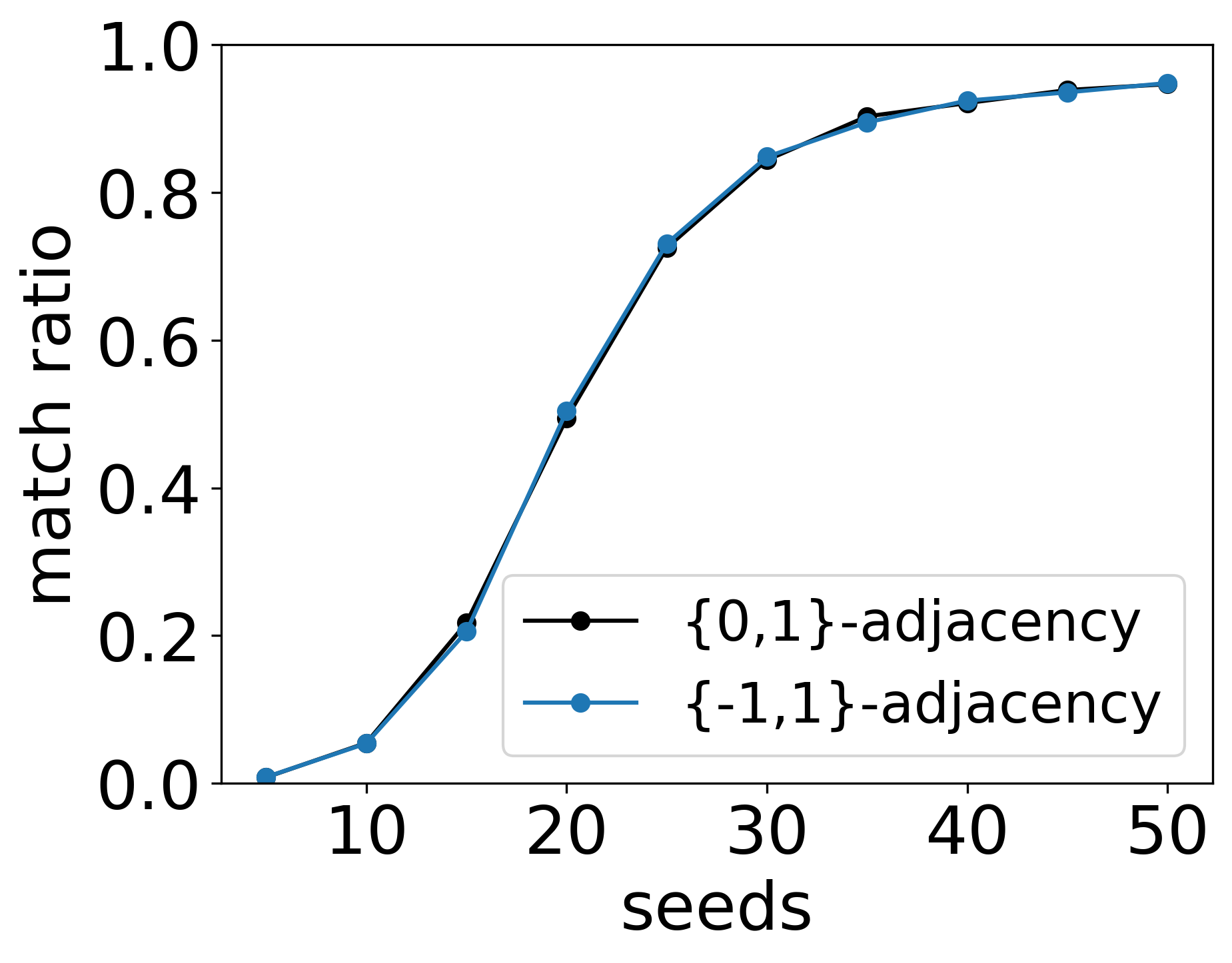}
    \caption*{$p = 0.5, \varrho = 0.7$}
    \end{subfigure}
    \caption{Comparing SGM when using two different forms of adjacency matrices.}
    \label{fig:SGM_adj}
\end{figure}

\newpage

\end{document}